\documentclass{nyj}
\usepackage{amssymb}
\usepackage{mathrsfs}
\usepackage[all]{xy}

\newtheorem{thm}{Theorem}[section]
\newtheorem{lem}[thm]{Lemma}
\newtheorem{lemma}[thm]{Lemma}

\newtheorem{prop}[thm]{Proposition}

\theoremstyle{definition}
\newtheorem{defn}[thm]{Definition}

\theoremstyle{remark}  
\newtheorem{rem}[thm]{Remark}
\newtheorem{remark}[thm]{Remark}

\numberwithin{equation}{section}

\renewcommand{\labelenumi}{(\roman{enumi})}

\newcommand{\C}{\mathbb C}

\newcommand{\KK}{\mathcal K}
\newcommand{\HH}{\mathcal H}
\newcommand{\GG}{\mathcal G}
\newcommand{\EE}{\mathcal E}
\renewcommand\SS{\mathcal{S}}

\newcommand\VVV{\mathscr{V}}

\renewcommand{\a}{\alpha}
\renewcommand{\b}{\beta}
\renewcommand{\d}{\delta}
\renewcommand{\l}{\lambda}
\renewcommand{\r}{\rho}
\newcommand{\p}{\phi}

\newcommand{\m}{\mu}

\newcommand{\s}{\sigma}

\renewcommand{\t}{\theta}
\renewcommand{\epsilon}{\varepsilon}
\renewcommand{\o}{\omega}
\newcommand{\z}{\zeta}

\newcommand{\G}{\Gamma}
\renewcommand{\P}{\Phi}
\renewcommand{\O}{\Omega}
\renewcommand{\L}{\Lambda}
\newcommand{\D}{\Delta}

\DeclareMathOperator{\aut}{Aut}

\DeclareMathOperator{\ad}{Ad}
\DeclareMathOperator{\supp}{supp}
\DeclareMathOperator{\lin}{Lin}

\DeclareMathOperator*{\spn}{span}
\DeclareMathOperator*{\clspn}{\overline{\spn}}

\newcommand{\id}{\text{\textup{id}}}

\newcommand{\inv}{^{-1}}

\renewcommand{\bar}{\overline}
\newcommand{\rt}{\textup{rt}}
\newcommand{\lt}{\textup{lt}}

\newcommand\bundlefont[1]{\mathscr{#1}}
\renewcommand\AA{\bundlefont{A}}
\newcommand\BB{\bundlefont{B}}
\def\sa_#1(#2,#3){\Gamma_{#1}(#2;#3)}
\def\cs(#1,#2){C^{*}(#1,#2)}
\newcommand\tA{\tilde A}
\newcommand\set[1]{\{\,#1\,\}}
\newcommand\X{\mathsf{X}}
\def\ipr(#1|#2){(#1 \mid  #2)_{\rho}}
\newcommand\tensor{\otimes}
\newcommand\tz{\t^{\AA}_{0}}
\newcommand\piz{\pi^{\AA}_{0}}
\newcommand\go{\GG^{(0)}}

\newcommand\gltg{G\times_{\lt}G}
\newcommand\altg{\AA\times_{\lt}G}
\newcommand\betart{\id_{G}\times\rt}
\newcommand\alphart{\id_{\AA}\times\rt}
\newcommand\alphabarrt{(\id_{\AA}\times\rt)\bar{\phantom{t}}}
\newcommand\csga{\cs(G,\AA)}
\newcommand\gltgg{\mathcal{S}}
\newcommand\altgg{\AA\times_{\lt}G\times_{\alphart}G}
\newcommand\btaa{\BB\times_{\alpha}G}

\newcommand\half{\frac12}
\newcommand\pan{\rho^{\AA}_{0}}
\newcommand\pg{\rho_{\GG}}
\newcommand\pa{\rho_{\AA}}
\newcommand\csfd{\cs(G,\AA)}
\newcommand\csmw{C^{*}_{\operatorname{Gr}}(G,\AA)}
\newcommand\safd{\sa_{c}(G,\AA)}
\newcommand\samw{\Gamma_{c}^{\operatorname{Gr}}(G,\AA)}
\let\phi\varphi 
\def\<{\langle}
\def\>{\rangle}
\let\ipscriptstyle=\scriptscriptstyle
\def\lipsqueeze{{\mskip -3.0mu}}
\def\ripsqueeze{{\mskip -3.0mu}}
\def\ipcomma{\nobreak\mathrel{,}\nobreak}
\newbox\ipstrutbox
\setbox\ipstrutbox=\hbox{\vrule height8.5pt
width 0pt}
\def\ipstrut{\copy\ipstrutbox}
\def\lip#1<#2,#3>{\mathopen{\relax_{\ipstrut\ipscriptstyle{
#1}}\lipsqueeze
\langle} #2\ipcomma #3 \rangle}
\def\blip#1<#2,#3>{\mathopen{\relax_{\ipstrut
\ipscriptstyle{ #1}}\lipsqueeze\bigl\langle} #2\ipcomma #3 \bigr\rangle}
\def\rip#1<#2,#3>{\langle #2\ipcomma #3
\rangle_{\ripsqueeze\ipstrut\ipscriptstyle{#1}}}
\def\brip#1<#2,#3>{\bigl\langle #2\ipcomma #3
\bigr\rangle_{\ripsqueeze\ipstrut\ipscriptstyle{#1}}}
\def\angsqueeze{\mskip -6mu}
\def\smangsqueeze{\mskip -3.7mu}
\def\trip#1<#2,#3>{\langle\smangsqueeze\langle #2\ipcomma #3
\rangle\smangsqueeze\rangle_{\ripsqueeze\ipstrut\ipscriptstyle{#1}}}
\def\btrip#1<#2,#3>{\bigl\langle\angsqueeze\bigl\langle #2\ipcomma
#3
\bigr\rangle
\angsqueeze\bigr\rangle_{\ripsqueeze\ipstrut\ipscriptstyle{#1}}}
\def\tlip#1<#2,#3>{\mathopen{\relax_{\ipstrut\ipscriptstyle{
#1}}\lipsqueeze \langle\smangsqueeze\langle} #2\ipcomma #3
\rangle\smangsqueeze\rangle}
\def\btlip#1<#2,#3>{\mathopen{\relax_{\ipstrut\ipscriptstyle{
#1}}\lipsqueeze
\bigl\langle\angsqueeze\bigl\langle} #2\ipcomma #3
\bigr\rangle\angsqueeze\bigr\rangle}

\def\ip(#1|#2){(#1\mid #2)}
\def\bip(#1|#2){\bigl(#1 \mid #2\bigr)}
\def\Bip(#1|#2){\Bigl( #1 \bigm| #2 \Bigr)}

\renewcommand\MR[1]{\relax}

\allowdisplaybreaks

\begin{document}

\title{Coactions and Fell bundles}

\author[Kaliszewski]{S. Kaliszewski}
\address{Department of Mathematics and Statistics
\\Arizona State University
\\Tempe, Arizona 85287}
\email{kaliszewski@asu.edu}

\author[Muhly]{Paul S. Muhly}
\address{Department of Mathematics
\\The University of Iowa
\\Iowa City, IA 52242}
\email{pmuhly@math.uiowa.edu}

\author[Quigg]{John Quigg}
\address{Department of Mathematics and Statistics
\\Arizona State University
\\Tempe, Arizona 85287}
\email{quigg@asu.edu}

\author[Williams]{Dana P. Williams}
\address{Department of Mathematics
\\Dartmouth College
\\Hanover, NH 03755}
\email{dana.williams@dartmouth.edu}



\begin{abstract}
We show that if $\AA$ is a Fell bundle over a locally compact group $G$, 
then there is a natural coaction $\delta$ of $G$ on the Fell-bundle
$C^*$-algebra $\csga$ 
such that if $\hat{\delta}$ is the dual action of $G$ on the crossed
product $\csga \rtimes_{\delta} G$, 
then the full crossed product $(\csga
\rtimes_{\delta}G)\rtimes_{\hat{\delta}}G$ 
is canonically isomorphic to $\csga\otimes\KK(L^2(G))$.
Hence the coaction $\delta$ is maximal.
\end{abstract}

\thanks{This research was partially funded by the Edward Shapiro fund
  at Dartmouth College.}

\thanks{Date: 15 September 2009}

\maketitle
\tableofcontents

\section*{Introduction}
\label{intro}

The theorem announced in the abstract, which we prove as Theorem
\ref{Phi-iso}, is part of a larger program that is inspired by the
realization, which only recently has come into focus, that Fell
bundles over groups and, more generally, Fell bundles over groupoids,
provide a natural setting for a broad range of imprimitivity theorems
and equivalence theorems for $C^*$-dynamical systems, especially
theorems involving nonabelian duality.  The present paper is a first
step in this larger program.

Very roughly, a Fell bundle $\AA$ over a locally compact group $G$ is
a bundle over $G$ such that the fibre $A_e$ over the identity $e$ of
$G$ is a $C^*$-algebra and such that the fibre $A_s$ over each $s\in G$ is
an $A_{e}$\,--\,$A_{e}$-imprimitivity bimodule with
the property that $A_s\otimes_{A_e} A_t$ is isomorphic to $A_{st}$ in
such a way that tensoring gives an associative multiplication on
$\AA$.\footnote{We follow the convention that the total space of a
  Banach bundle is represented in a script font, while the fibres are
  written in Roman font.  Thus if $p:\AA\to X$ is a bundle over a
  space $X$, then we'll write $A_x$ for the fibre $p^{-1}(x)$ viewed
  as a Banach space.}  The space of continuous, compactly supported
cross sections of $\AA$, denoted $\Gamma_c(G;\AA)$, carries a natural
convolution-like product under which it forms a $*$-algebra. A certain
completion of this algebra is a $C^*$-algebra, denoted $\csga$.  One
can profitably think of $\csga$ as a generalized crossed product of
$A_e$ by $G$.  Indeed, if $G$ acts on a $C^*$-algebra $B$ via a
continuous homomorphism $\alpha : G \to \aut (B)$, and if $\AA$ is
defined to be $B\times G$, with product defined by the equation
$(a,s)(b,t)=(a\alpha_{s}(b),st)$, then $\AA$ is a Fell bundle over $G$, called
the \emph{semidirect-product bundle} determined by the action,
and the $C^*$-crossed product $B\rtimes_\alpha G$ is isomorphic to the
bundle $C^*$-algebra $\csga$. This point was made by Fell in his first
works on the subject \cite{fel:ajm69,fel:mams69} and was one of the
reasons he began the theory of these bundles.  Importantly, not every
Fell bundle over a group $G$ is isomorphic to such a bundle \cite[\S\S
VIII.3.16, VIII.4.7]{fd2}.\footnote{We shall have more to say about
semidirect-product bundles in Sections \ref{Semidirect-product bundles}~and
\ref{action crossed product}.}

Coactions were introduced to give a generalization, for non-abelian
groups, of the Takai-Takesaki duality for crossed products by actions
of abelian groups on $C^*$-algebras.  Subsequently Katayama proved a
crossed-product duality theorem for coactions, specifically, if
$\delta$ is a coaction of a group $G$ on a $C^*$-algebra $A$, then
there is a dual action $\hat\delta$ of $G$ on the crossed product
$A\rtimes_\delta G$ such that the \emph{reduced} crossed product
$(A\rtimes_\delta G)\rtimes_{\hat\delta,r}G$ is isomorphic to
$A\otimes \KK(L^2(G))$.  Katayama used what are now known as
\emph{reduced} coactions, which involve the reduced group
$C^*$-algebra $C^*_r(G)$.  For more information on crossed-product
duality, see \cite[Appendix~A]{enchilada}.

The use of the term ``crossed product'' both in the context of group
actions and in the context of coactions may seem confusing, initially.
However, in practice, it is easy to distinguish between the two.

Raeburn introduced \emph{full} coactions, which involve the full group
$C^*$-algebra $C^*(G)$, to take advantage of universal properties.
For such coactions, there is always a canonical surjection
\[
\Phi: A\rtimes_\delta G\rtimes_{\hat\delta}G\to A\otimes \KK(L^2(G)),
\]
and the question naturally arose, when is $\Phi$ in fact an
isomorphism?  When this is the case, \emph{full crossed-product
  duality} is said to hold, and the coaction $\delta$ is said to be
\emph{maximal}. For example, the dual coaction on a full
crossed product by an action is always maximal
\cite[Proposition~3.4]{ekq}.

Since Fell bundle $C^*$-algebras are generalizations of crossed
products by actions, it is natural to ask whether there exists a
coaction $\d$ of $G$ on $C^*(G,\AA)$, and if so, whether $\delta$ is
maximal.  In the present paper, we settle these questions
affirmatively.

The existence of a coaction on $\csfd$ was briefly presented in
\cite{lprs} for the case of reduced coactions.  In \cite{Q96} (see
also \cite{ngdiscrete}), the third author showed that when the group
$G$ is discrete there is in fact a bijective correspondence between
Fell bundles over $G$ and coactions of $G$ on $C^*$-algebras.
Further, in \cite{eq:full} the third author and Echterhoff observed
that given a Fell bundle $\AA$ over a discrete group $G$, there is a
natural coaction $\delta$ of $G$ on $\csga$ and the crossed product
$\csga \rtimes_{\delta} G$ is naturally isomorphic to the
$C^*$-algebra of a Fell bundle $\AA \times_\lt G$ over the {discrete}
groupoid $G\times_{\lt}G$ obtained by letting $G$ act on itself by
left translation.  This observation, coupled with the work of the
second and fourth authors on the theory of Fell bundles over groupoids
\cite{mw:fell} (which, in turn, was inspired, in part, by \cite{Q96}),
was the point of departure for the current project.

Indeed, although groupoids do not appear explicitly in the statement
of our main theorem, Fell bundles over groupoids are crucial in the
techniques we develop for the proof.  We rely heavily on
\cite{mw:fell} for the theory and basic results concerning Fell
bundles over groupoids.  In particular, we make free use of the
Disintegration Theorem for Fell bundles \cite[Theorem~4.13]{mw:fell}
which is a generalization of Renault's Disintegration Theorem for
groupoids \cite[Proposition~4.2]{ren:representation}. (See
\cite[\S7]{muhwil:nyjm08} for more discussion and references on
Renault's Theorem.)

The plan for our proof of Theorem \ref{Phi-iso} is as follows: The
initial two sections are preparatory. Section \ref{prelim} establishes
notation and collects some results that will be used in the sequel.
Section \ref{product bundles} addresses some fine points regarding the
problem of ``promoting'' a Fell bundle over a group to a Fell bundle
over the product of the group with itself.  The first real step in our
analysis is taken in Section \ref{exist}. There we prove in
Proposition \ref{exists} that if $\AA$ is a Fell bundle over a locally
compact group $G$, then there is a natural coaction $\delta$ of $G$ on
$\csga$ analogous to the dual coaction on a crossed product.

We note in passing that in \cite{exelng}, Exel and Ng prove a result
that is similar to our Proposition \ref{exists}. However, their
setting is somewhat different from ours in that it uses an older and
no-longer-used definition of ``full coaction'' that was advanced by
Raeburn in \cite{rae:representation}.  Also, their proof is different
in certain important respects.  So, to keep this note self-contained
we present full details.

The second substantial step taken in our analysis is Theorem
\ref{coaction crossed product isomorphism}, which asserts that there
is a natural isomorphism $\theta$ from the crossed product $\csga
\rtimes_\delta G$ to the $C^*$-algebra $C^*(G\times_{\lt}G,\AA
\times_\lt G)$ of the Fell bundle $\AA \times_\lt G$ over the
transformation groupoid $G\times_{\lt}G$.  As we mentioned above, this
isomorphism theorem was inspired by \cite{eq:full}. Section
\ref{transformation bundles} provides the necessary prerequisites for
the formulation and proof of Theorem \ref{coaction crossed product
  isomorphism}.

The third major step is Theorem~\ref{action crossed product
  isomorphism}, which establishes, in the general context of a Fell bundle $\BB$
over a groupoid $\GG$,
an isomorphism between the $C^*$-algebra of a
semidirect-product bundle $\BB\times_\alpha G$ (the theory of which is
developed in Section~\ref{Semidirect-product bundles}) and the crossed
product of $C^*(\GG,\BB)$ by a corresponding action of~$G$.

The remainder of the argument occupies Section~\ref{canonical}.
There, we show that the isomorphism $\theta$ established in Theorem
\ref{coaction crossed product isomorphism} is equivariant for the dual
action $\hat{\delta}$ of~$G$ on $\csga \times_\delta G$ and a natural
action of~$G$ on $C^*(G\times_{\lt}G,\AA \times_\lt G)$.  Using this,
$\theta$ is promoted to an isomorphism between the two crossed
products.  We then apply the result of Section~\ref{action crossed
  product} to this natural action to see that the crossed product can
be realized as the $C^*$-algebra of a certain semidirect-product
bundle; this bundle turns out to be isomorphic to one whose
$C^*$-algebra is easily recognized as $C^*(G,\AA)\otimes\KK(L^2(G))$.
Finally, we show that these isomorphisms combine to give the canonical
surjection~$\Phi$, and this completes our proof of Theorem
\ref{Phi-iso}.

\section{Preliminaries}
\label{prelim}

If $A$ is a $C^{*}$-algebra, then its maximal unitization $M(A)$
(\cite[Definition~2.46]{tfb}) is called the \emph{multiplier algebra}
of $A$.  Traditionally, $M(A)$ is realized as the collection of double
centralizers.  Here we adopt the approach taken in \cite{tfb},
regarding $M(A)$ as the algebra $\mathcal{L}(A)$ of bounded
adjointable operators on $A$ viewed as a right-Hilbert module over
itself.  (That any two maximal unitizations are naturally isomorphic
is guaranteed by \cite[Theorem~2.47]{tfb}.)  As usual, we let $\tA$ be
the $C^{*}$-subalgebra of $M(A)$ generated by $A$ and $1_{M(A)}$.  (Thus
$\tA=A$ if $A$ is unital,
and $\tA$ is $A$ with an identity adjoined
otherwise.)  We use minimal tensor products of $C^*$-algebras
throughout.

Let $G$ be a locally compact group.  We use $u: G\to M(C^*(G))$ to
denote the canonical embedding, although sometimes we will simply
identify $s\in G$ with its image $u(s)\in M(C^*(G))$.  Similarly, we
will usually not distinguish between a strictly continuous unitary
homomorphism of $G$ and its unique nondegenerate extension to
$C^*(G)$.  As a general reference for group actions we use
\cite{danacrossed}, and for coactions we refer to
\cite[Appendix~A]{enchilada}.

\subsection{Group Actions}
An \emph{action} of $G$ on a $C^*$-algebra $A$ is a homomorphism
$\alpha: G\to \aut A$ such that the map $s\mapsto \alpha_s(a)$ is norm
continuous from $G$ to $A$ for each $a\in A$.  A \emph{covariant
  representation} of $(A,G,\alpha)$ on a Hilbert space $\HH$ is a pair
$(\pi,U)$, where $\pi: A\to B(\HH)$ is a nondegenerate representation
and $U: G\to B(\HH)$ is a strongly continuous unitary representation,
which satisfies the \emph{covariance condition}
\begin{equation}
  \label{variant}
  \pi(\alpha_s(a)) = U_s \pi(a) U_s^*
  \quad\text{for $a\in A$ and $ s\in G$.}
\end{equation}
More generally, for any $C^*$-algebra $B$, a \emph{covariant
  homomorphism} of $(A,G,\alpha)$ into $M(B)$ is a pair $(\pi,U)$,
where $\pi:A\to M(B)$ is a nondegenerate homomorphism and $U: G\to
M(B)$ is a strictly continuous unitary homomorphism, which
satisfies~\eqref{variant}.

A \emph{crossed product} for $(A,G,\alpha)$ is a $C^*$-algebra
$A\rtimes_\alpha G$, together with a covariant homomorphism
$(i_A,i_G)$ of $(A,G,\alpha)$ into $M(A\rtimes_\alpha G)$ which is
universal in the sense that for any covariant homomorphism $(\pi,U)$
of $(A,G,\a)$ into $M(B)$ there is a unique nondegenerate homomorphism
$\pi\rtimes U:A\rtimes_\a G\to M(B)$, called the \emph{integrated
  form} of $(\pi,U)$, such that
\[
\pi=(\pi\rtimes U)\circ i_A\quad\text{and}\quad U=(\pi\rtimes U)\circ
i_G.
\]
The crossed product is generated by the universal covariant
homomorphism in the sense that
\[
A\rtimes_\a G = \clspn\set{ i_A(a) i_G(f) : \text{$a\in A$ and $ f\in
    C_c(G)$}}.
\]
The space $C_c(G,A)$ of compactly supported continuous functions from
$G$ into $A$ is a $*$-algebra with (convolution) multiplication and
involution given by
\[
(f*g)(s) = \int_G f(t) \alpha_t(g(t\inv s)) \,dt \quad\text{and}\quad
f^*(s) = f(s\inv)^*\D(s)\inv,
\]
where $\D$ denotes the modular function of $G$.  The algebra
$C_c(G,A)$ embeds as a dense $*$-subalgebra of $A\rtimes_\alpha G$ via
the map
\[
f\mapsto \int_G i_A(f(s)) i_G(s)\,ds,
\]
so that if $(\pi,U)$ is a covariant homomorphism of $(A,G,\alpha)$,
then
\[
\pi\rtimes U(f) = \int_G \pi(f(s)) U(s)\,ds.
\]

\subsection{Coactions}
A \emph{coaction} of $G$ on a $C^*$-algebra $A$ is a nondegenerate
injective homomorphism $\d:A\to M(A\otimes C^*(G))$ which satisfies
the \emph{coaction identity}
\begin{equation}
  \label{comodule}
  (\d\otimes \id_G)\circ \d=(\id\otimes \d_G)\circ \d,
\end{equation}
and which is \emph{nondegenerate as a coaction} in the sense that
\begin{equation}
  \label{nondegenerate}
  \clspn\{\d(A)(1\otimes C^*(G))\}=A\otimes C^*(G).
\end{equation}
Here $\d_G:C^*(G)\to M(C^*(G)\otimes C^*(G))$ is the homomorphism
determined by the unitary homomorphism of $G$ given by $s\mapsto
u(s)\otimes u(s)$.  Note that condition~\eqref{nondegenerate} implies
nondegeneracy of $\d$ as a map into $M(A\otimes C^*(G))$.

A \emph{covariant representation} of $(A,G,\d)$ on a Hilbert space
$\HH$ is a pair $(\pi,\m)$, where $\pi: A\to B(\HH)$ and $\m:
C_0(G)\to B(\HH)$ are nondegenerate representations which satisfy the
\emph{covariance condition}
\begin{equation}
  \label{covariant}
  \ad (\m\otimes\id)(w_G)(\pi(a)\otimes 1)=(\pi\otimes\id)(\d(a))
  \quad\text{for $a\in A$.}
\end{equation}
Here $w_G$ is the element of $M(C_0(G)\otimes C^*(G))$ which
corresponds to the canonical embedding $u: G\to M(C^*(G))$ under the
natural isomorphism of $M(C_0(G)\otimes C^*(G))$ with the strictly
continuous bounded maps from $G$ to $M(C^*(G))$.  More generally, for
any $C^*$-algebra $B$, a \emph{covariant homomorphism} of $(A,G,\d)$
into $M(B)$ is a pair $(\pi,\m)$, where $\pi:A\to M(B)$ and
$\m:C_0(G)\to M(B)$ are nondegenerate homomorphisms satisfying
\eqref{covariant}.

A \emph{crossed product} for $(A,G,\d)$ is a $C^*$-algebra
$A\rtimes_\d G$, together with a covariant homomorphism $(j_A,j_G)$ of
$(A,G,\d)$ into $M(A\rtimes_\d G)$ which is universal in the sense
that for any covariant homomorphism $(\pi,\m)$ of $(A,G,\d)$ into
$M(B)$ there is a unique nondegenerate homomorphism $\pi\rtimes
\m:A\rtimes_\d G\to M(B)$, called the \emph{integrated form} of
$(\pi,\mu)$, such that
\[
\pi=(\pi\rtimes \m)\circ j_A\quad\text{and}\quad \m=(\pi\rtimes
\m)\circ j_G.
\]
The crossed product is generated by the universal covariant
homomorphism in the sense that
\[
A\rtimes_\delta G = \clspn\set{ j_A(a) j_G(f) : \text{$a\in A$ and $
    f\in C_0(G)$}}.
\]
The \emph{dual action} of $G$ on $A\rtimes_\delta G$ is the
homomorphism $\hat\delta: G\to \aut(A\rtimes_\d G)$ given on
generators by
\[
\hat\delta_s( j_A(a) j_G(f) ) = j_A(a) j_G(\rt_s(f)),
\]
where $\rt$ denotes the action of $G$ on $C_0(G)$ by right
translation: $\rt_s(f)(t) = f(ts)$.

Given a representation $\pi$ of $A$ on a Hilbert space $\HH$, the
associated \emph{regular representation} $\L$ of $A\rtimes_\delta G$
on $\HH\otimes L^2(G)$ is the integrated form
\[
\L=\bigl((\pi\otimes\l)\circ \d\bigr)\rtimes (1\otimes M),
\]
where $\lambda$ is the left regular representation of $G$ on $L^2(G)$
and $M$ is the representation of $C_0(G)$ on $L^2(G)$ by
multiplication: $(M_f\xi)(s) = f(s)\xi(s)$.  When $\pi$ is faithful,
the associated regular representation is always faithful
\cite[Remark~A.43(3)]{enchilada}, and thus gives an isomorphism
between $A\rtimes_\delta G$ and the concrete $C^*$-algebra
\[
\L(A\rtimes_\delta G) = \clspn\set{ (\pi\otimes\l)\circ \d(a)(1\otimes
  M_f) : \text{$a\in A$ and $ f\in C_0(G)$} }.
\]

The \emph{canonical surjection} associated to $\delta$ is the map
\[
\Phi = \bigl( (\id\otimes\lambda)\circ\delta \rtimes (1\otimes
M)\bigr) \rtimes (1\otimes\rho): A\rtimes_\delta G\rtimes_{\hat\delta}
G \to A\otimes \KK(L^2(G)),
\]
where $\rho$ is the right regular representation of $G$ on $L^2(G)$.
(It almost goes without saying that, by convention,
$(\lambda_s\xi)(t) = \xi(s\inv t)$ and
$(\rho_s\xi)(t) = \xi(ts)\D(s)^{1/2}$.)
On the generators, $\Phi$ is given by
\[
\Phi\bigl( i_{A\rtimes_\delta G}(j_A(a)j_G(f)) i_G(g) \bigr) =
(\id\otimes\lambda)\circ\delta(a)\bigl( 1\otimes M_f\rho(g)\bigr)
\]
for $a\in A$, $f\in C_0(G)$, and $g\in C^*(G)$.  The coaction $\delta$
is \emph{maximal} if the canonical surjection $\Phi$ associated to
$\delta$ is injective; thus the maximal coactions are precisely those
coactions for which full crossed-product duality holds in the sense
that $\Phi$ is an isomorphism of $A\rtimes_\delta
G\rtimes_{\hat\delta} G$ onto $A\otimes \KK(L^2(G))$.

Some of our coaction calculations will involve the
\emph{Fourier-Stieltjes} algebra $B(G)$.  (see \cite[\S\S
A.4--A.5]{enchilada} for brief survey or \cite{eym} for a more
detailed treatment).  In simple terms the Fourier-Stieltjes algebra
$B(G)$ is a space of bounded continuous functions on $G$ which can be
identified with the dual space $C^*(G)^*$ via the formula
\[
f(g) = \int_G f(s) g(s) \,ds \quad\text{for $f\in B(G)$ and $ g\in
  C_c(G)\subseteq C^*(G)$.}
\]
By \cite[Propositions~3.4 and~3.7]{eym}, 
the intersection $B(G)\cap C_0(G)$ is norm dense in~$C_0(G)$.
For $f\in B(G)$, the \emph{slice map} $\id_{A}\otimes f: A\otimes
C^*(G)\to A$ determined by
\[
(\id_{A}\otimes f)(a\otimes b)=af(b)\quad\text{for $a\in A$ and $b\in
  C^*(G)$}
\]
extends uniquely to a strictly continuous linear map $\id_{A}\otimes
f: M(A\otimes C^*(G))\to M(A)$, and moreover such slice maps separate
the points of $M(A\otimes C^*(G))$ (\cite[Lemma~A.30]{enchilada}).

\subsection{Fell Bundles}
A Fell bundle over a groupoid is a natural generalization of Fell's
$C^{*}$-algebraic bundles over groups treated in detail in
\cite[Chap.~~VIII]{fd2} and discussed briefly in the introduction.  We
will refer to \cite{mw:fell} for the particulars of Fell bundles over
groupoids.  Generally speaking, a \emph{Fell bundle} $p:\BB\to\GG$ is
a upper semicontinuous
 Banach bundle over a locally compact Hausdorff groupoid
$\GG$ satisfying the axioms laid out in
\cite[Definition~1.1]{mw:fell}.\footnote{There are a number of
  equivalent definitions of Fell bundles over groupoids in the
  literature starting with Yamagami's original in
  \cite[Definition~1.1]{yam:symmetric}, as well as
  \cite[Definition~6]{muh:cm01} and \cite[Definition~2.1]{dkr:ms08}.}
It was observed in \cite[Lemma~3.30]{bmz} that the underlying
    Banach bundle of an upper semicontinuous Fell bundle over a
    \emph{group} is necessarily continuous.  (The authors of
    \cite{bmz} attribute this observation to Exel.) Since all the Fell
    bundles in this work originate from Fell bundles over groups, they
    will necessarily be built on continuous Banach
    bundles.\footnote{An exception is that in sections
      \ref{Semidirect-product bundles}~and 
      \ref{action crossed product} we work with general Fell bundles
      over groupoids, and there it is not necessary to assume that the
      underlying Banach bundles are continuous.}
We will
assume all the Fell bundles here are \emph{separable} in that $\GG$ is
second countable and the Banach space $\sa_{0}(\GG,\BB)$ of sections
is separable. (This hypothesis is not only a sign of good taste, but
it will also ensure that the results of \cite{mw:fell} apply.)

We are only interested in groupoids $\GG$ with a continuous Haar
system $\set{\lambda^{u}}_{u\in\go}$.  Then the set $\sa_{c}(\GG,\BB)$ of
continuous compactly supported sections 
of $\BB$ has
the structure of a $*$-algebra:
\begin{equation*}
  f*g(x):=\int_{\GG}f(y)g(y^{-1}x)\,d\lambda^{r(x)}(y)\quad \text{and}
  \quad f^{*}(x):=f(x^{-1})^{*}. 
\end{equation*}
Then we can define a norm, $\|\cdot\|_{I}$, on $\sa_{c}(\GG,\BB)$ via
\begin{equation*}
  \|f\|_{I}=\max\Bigl\{ \sup_{u\in\go}\int_{\GG}
  \|f(x)\|\,d\lambda^{u}(x),
  \sup_{u\in\go}\int_{G}\|f(x)\|\,\lambda_{u}(x) \Bigr\}.
\end{equation*}
If $\HH$ is a Hilbert space, then a $*$-homomorphism
$L:\sa_{c}(\GG,\BB)\to B(\HH)$ is called
\emph{$\|\cdot\|_{I}$-decreasing} if $\|L(f)\|\le\|f\|_{I}$ for all
$f$.  We say that $L$ is a
  $\|\cdot\|_{I}$-decreasing \emph{representation} if it is also
  \emph{nondegenerate} in the sense that
\begin{equation*}
  \overline{\operatorname{span}}\set{L(f)\xi:\text{$f\in\sa_{c}(\GG,\BB)$
      and $\xi\in\HH$}}=\HH.
\end{equation*}
Then, by definition, the \emph{universal norm} on $\sa_{c}(\GG,\BB)$
is
\begin{equation*}
  \|f\|:=\sup\set{\|L(f)\|:\text{$L$ is a $\|\cdot\|_{I}$-decreasing
      representation of $\sa_{c}(\GG,\BB)$}}.
\end{equation*}
The completion $\overline{\bigl(\sa_{c}(\GG,\BB),\|\cdot\|\bigr)}$ is
the $C^{*}$-algebra $\cs(\GG,\BB)$ of the Fell bundle
$p:\BB\to\GG$.\footnote{It might be helpful to look over the examples
  in \cite[\S2]{mw:fell} at this point.}

More generally, a nondegenerate $*$-homomorphism $L:\G_c(\GG,\BB)\to
B(\HH)$ is called simply a \emph{representation} if $L$ is continuous
when $\sa_{c}(\GG,\BB)$ is equipped with the inductive limit topology
and $B(\HH)$ is given the weak operator topology.  It is a nontrivial
result --- a consequence of the Disintegration Theorem
(\cite[Theorem~4.13]{mw:fell}) --- that every representation of
$\sa_{c}(\GG,\BB)$ is $\|\cdot\|_{I}$-decreasing.  Since
$\|\cdot\|_{I}$-decreasing representations are clearly
representations, we see that
\begin{equation*}
  \|f\|=\sup\set{\|L(f)\|:\text{$L$ is a 
      representation of $\sa_{c}(\GG,\BB)$}}
\end{equation*}
(see \cite[Remark~4.14]{mw:fell}).

\begin{lemma}
  \label{lem-pull-back}
  Suppose that $p:\BB\to \GG$ is a Fell bundle over a locally compact
  groupoid $\GG$.  If $\HH$ is a locally compact groupoid and
  $\phi:\HH\to \GG$ is a continuous groupoid homomorphism, then the
  pull-back $q:\phi^{*}\BB\to \HH$ is a Fell bundle over $\HH$ with
  multiplication and involution given by
  \begin{equation*}
    (a,h)(b,t)=(ab,ht)\quad\text{and}\quad (a,h)^{*}=(a^{*},h^{-1}).
  \end{equation*}
\end{lemma}
\begin{proof}
  The proof is routine.  For example, $q:\phi^{*}\BB\to \HH$ is
  clearly a Banach bundle (see \cite[\S II.13.7]{fd1} where pull-backs
  are called retractions).  The fibre over $h$ is isomorphic to
  $B_{\phi(h)}$.  The Fell bundle structure from $\BB$ makes the
  latter into a $B_{r(\phi(h))}$\,--\,$B_{s(\phi(h))}$-imprimitivity
  bimodule.  Since the fibre over $s(h)$ is isomorphic to
  $B_{\phi(s(h))}$ and $\phi(s(h))=s(\phi(h))$, the rest is easy.
  (Note that when $\GG$ and $\HH$ are groups, this result is \cite[\S
  VIII.3.17]{fd2}.)
\end{proof}

\subsection{Fell Bundles over Groups}
However, to begin with, we are interested in a (separable, of course)
Fell bundle $p:\AA\to G$ where $G$ is a locally compact \emph{group}.
This case affords a number of simplifications, and also allows us to
avoid some of the overhead coming from \cite{mw:fell}.  Note that a
Fell bundle $p:\AA\to G$ over a group is what Fell and Doran call a
$C^{*}$-algebraic bundle over $G$ (see \cite[Definitions~VIII.16.2 and
VIII.3.1]{fd2}).  Since we ultimately treat Fell bundles over groups
as a special case of a Fell bundle over a groupoid, our axioms require
that $p:\AA\to G$ is \emph{saturated} in the sense that
$\clspn\{A_sA_t\}=A_{st}$ for all $s,t\in G$ (see \cite[\S
VIII.2.8]{fd2}).  We will often write $a_{s}$ for an element of
$A_{s}$; that is, $a_{s}\in \AA$ and $p(a_{s})=s$.

We do make one deviation from the groupoid treatment when building the
associated $C^{*}$-algebra, $\cs(G,\AA)$.  In order that we can easily
obtain the usual group $C^{*}$-algebra construction as well as the
usual crossed-product construction as special cases, it is convenient
to add the modular function, $\Delta$, on $G$ to the definition of the
involution on $\sa_{c}(G,\AA)$:
\begin{equation*}
  f^{*}(s)=\Delta(s)^{-1}f(s^{-1})^{*}
\end{equation*}
(see \cite[\S VIII.5.6]{fd2}).  Then the somewhat unsatisfactory
$\|\cdot\|_{I}$ reduces to the normal analog of the $L^{1}$-norm:
\begin{equation*}
  \|f\|_{1}:=\int_{G}\|f(s)\|\,ds,
\end{equation*}
and the universal norm on $\sa_{c}(G,\AA)$ is given as the supremum
over $\|\cdot\|_{1}$-decreasing representations. As we shall see
shortly (see Remark~\ref{rem-modular-differences}), the isomorphism
class of $\cs(G,\AA)$ is the same as that obtained using the definition of
the involution given for groupoids where no modular function is
available.

Assuming $p:\AA\to G$ is a Fell bundle over a group, a
$*$-homomorphism $\pi:\AA\to M(B)$ is just a map with the obvious
algebraic properties.  We call $\pi$ \emph{nondegenerate} if
\begin{equation*}
  \clspn\{\pi(A_{e})B\}=B.
\end{equation*}
The next lemma shows that $\AA$ comes with a canonical nondegenerate
strictly continuous embedding $\iota:\AA\to M(\cs(G,\AA))$.  Then
Lemma~\ref{bundle map} shows that the pair $(C^*(G,\AA),\iota)$ is in
fact universal for strictly continuous nondegenerate $*$-homomorphisms
of $\AA$ into multiplier algebras.

\begin{lem}\label{key}
  Let $p:\AA\to G$ be a separable Fell bundle over a locally compact
  group $G$.  There exists a strictly continuous nondegenerate
  $*$-homomorphism $\iota: \AA\to M(\cs(G,\AA))$ such that for $a_s\in
  A_s$ and $f\in \sa_c(G,\AA)$, we have $\iota(a_s)f\in \sa_c(G,\AA)$,
  with
  \begin{equation}\label{iota}
    (\iota(a_s)f)(t) = a_s f(s\inv t).
  \end{equation}
\end{lem}

\begin{proof}
  For each $a_s\in A_s$, \eqref{iota} clearly defines a linear map
  $\iota(a_s)$ of $\sa_c(G,\AA)$ into itself.  Here we will view
  $\sa_{c}(G,\AA)$ as a dense subspace of $\cs(G,\AA)$ viewed as a
  Hilbert module over itself.  Then the inner product $\rip< f, g>
  =f^{*}*g$ is $\sa_{c}(G,\AA)$-valued on $\sa_{c}(G,\AA)$.  
It is easy to check that
  $\iota(a_{s})\iota(a_{t})=\iota(a_{s}a_{t})$, and a straightforward
  computation shows that
  \begin{equation}
    \label{eq:1}
    \brip< \iota(a_{s})f,g> = \brip< f,\iota(a_{s}^{*})g>
  \end{equation}
  (a similar, but more involved computation is given in detail in the
  proof of Theorem~\ref{coaction crossed product isomorphism}).  Since
  $\|a_{s}\|^{2}1_{A_{e}}-a_{s}^{*}a_{s}\ge 0$ in $\tA_{e}$, there is
  a $b_{e}\in \tA_{e}$ such that
  $\|a_{s}\|^{2}1_{A_{e}}-a_{s}^{*}a_{s} = b_{e}^{*}b_{e}$. Then,
since~\eqref{iota} makes sense and $\iota$ is multiplicative
for elements of $\tA_{e}$,
and since~\eqref{eq:1} also holds for $b_e\in\tA_{e}$,
we see that
  \begin{align*}
    \|a_{s}\|^{2}\brip< f,f> - \brip< \iota(a_{s})f,\iota(a_{s})f > &=
    \brip<
    \iota(\|a_{s}\|^{2}1_{A_{e}} -a_{s}^{*}a_{s})f,f> \\
    &= \brip< \iota(b_{e})f,\iota(b_{e})f> \ge 0
  \end{align*}
for all $f\in \sa_c(G,\AA)$.  
  It follows that $\iota(a_{s})$ is bounded and extends to a bounded
  operator on $\cs(G,\AA)$ with adjoint $\iota(a_{s}^{*})$.  It is
  routine to verify that the resulting map $\iota:\AA\to
  M(\cs(G,\AA))$ is a $*$-homomorphism.

  To see that $\iota$ is nondegenerate, first note that $A_{s}$ is an
  $A_{e}$\,--\,$A_{e}$-imprimitivity bimodule.  Thus if
  $\set{a_{i}}_{i\in I}$ is an approximate identity in $A_{e}$, then
  $a_{i}a_{s}\to a_{s}$ for any $a_{s}\in A_{s}$.  Then a messy
  compactness argument similar to that given in the proof of
  Theorem~\ref{coaction crossed product isomorphism} shows that
  $\iota(a_{i})f\to f$ in the inductive limit topology on
  $\sa_{c}(G,\AA)$ for any $f\in \sa_{c}(G,\AA)$.  Since convergence
  in the inductive limit topology implies convergence in the
  $C^{*}$-norm, this establishes nondegeneracy.

  It only remains to prove strict continuity.  Our separability
  assumptions on $p: \AA\to G$ allow us to invoke
  \cite[Proposition~II.13.21]{fd1} to see that $\AA$ is second
  countable.  Thus, it suffices to show that if $\set{a_{s_n}}$ is a
  sequence in $\AA$ converging to $a_s$, then $\iota(a_{s_n})\to
  \iota(a_s)$ strictly.

  The convergent sequence $\set{a_{s_n}}$ must lie in a norm-bounded
  subset of $\AA$, so the image $(\iota(a_{s_n}))$ is a bounded
  sequence in $M(C^*(G,\AA))$ (because $\|\iota(a_s)\|\leq \|a_s\|$).
  Thus, it suffices to show that $\iota(a_{s_n})\to \iota(a_s)$
  $*$-strongly; and since $a_{s_n}^*\to a_s^*$ and $\iota$ is
  $*$-preserving, it suffices to show strong convergence.  Finally,
  since $\set{\iota(a_{s_n})}$ is bounded, it suffices to show that
  $\iota(a_{s_n})f\to \iota(a_s)f$ in the inductive limit topology,
  for each $f\in \sa_c(G,\AA)$.

  Suppose not; so there is $f\in\sa_c(G,\AA)$ such that
  $\iota(a_{s_n})f$ does not converge to $\iota(a_s)f$ in the
  inductive limit topology.  Note that since $s_n\to s$ in $G$, we can
  find a compact set $K\subseteq G$ such that the supports of
  $\iota(a_s)f$ and all the $\iota(a_{s_n})f$ are contained in $K$, so
  it must be that the convergence is not uniform on $K$.  So, passing
  to a subsequence and relabeling, we can find $\epsilon>0$ and
  $t_n\to t$ in $K$ such that for all $n$,
  \[
  \| \iota(a_{s_n})f(t_n) - \iota(a_s)f(t_n) \| \geq \epsilon.
  \]
  But by joint continuity of multiplication in $\AA$, we have
  \[
  \iota(a_{s_n})f(t) = a_{s_n}f(s_n\inv t) \to a_s f(s\inv t) =
  \iota(a_s)f(t)
  \]
  in $\AA$.  Since this implies that the norm of the difference goes
  to zero, we have a contradiction.
\end{proof}

\begin{lem}
  \label{bundle map}
  Let $p:\AA\to G$ be as in Lemma~\ref{key}.  If $B$ is a
  $C^*$-algebra and $\pi_0:\AA\to M(B)$ is a strictly continuous
  nondegenerate $*$-homomorphism, then there is a unique nondegenerate
  homomorphism $\pi:\cs(G,\AA)\to M(B)$, called the \emph{integrated
    form} of $\pi_0$, such that $\pi\circ\iota = \pi_0$.  Moreover,
  \begin{equation}\label{int}
    \pi(f)=\int_G \pi_0(f(s))\,ds
    \quad\text{for $f\in \sa_c(G,\AA)$.}
  \end{equation}

  Conversely, every nondegenerate $*$-homomorphism of $C^*(G,\AA)$ is
  the integrated form of some such $\pi_0$.
\end{lem}

\begin{remark}
  Note that the integral in~\eqref{int} makes sense since $\pi_0\circ
  f$ is strictly continuous so that we can apply, for example,
  \cite[Lemma~C.11]{tfb}.
\end{remark}

\begin{proof}
  It is straightforward to check that \eqref{int} defines a
  $*$-homomorphism $\pi:\sa_c(G,\AA)\to M(B)$.

  To see that $\pi$ is nondegenerate, we need to see that
  \begin{equation*}
    \spn\set{\pi(f)b:\text{$f\in\sa_{c}(G,\AA)$ and $b\in B$}}
  \end{equation*}
  is dense in $B$.  To this end, fix $a\in A_{e}$ and choose
  $f\in\sa_{c}(G,\AA)$ such that $f(e)=a$.  Let $\set{\phi_{k}}$ be a
  sequence in $C_{c}^{+}(G)$ with integral one whose supports shrink
  to the identity.  Let $f_{k}(s)=\phi_{k}(s)f(s)$.  Then it is not
  hard to see that $\pi(f_{k})b\to \pi_{0}(a)b$.  Therefore, the
  nondegeneracy of $\pi$ follows from that of $\pi_{0}$.

  If $L:B\to B(\mathcal{H})$ is a faithful representation, then
  $L\circ \pi$ is a $\|\cdot\|_{1}$-decreasing representation of
  $\sa_{c}(G,\AA)$.  By the definition of the universal norm,
  \begin{equation*}
    \|L\circ \pi(f)\|\le\|f\|.
  \end{equation*}
  Since the extension of $L$ to $M(B)$ is isometric, $\|\pi(f)\|\le
  \|f\|$.  Therefore, $\pi$ extends to $\cs(G,\AA)$.

  To prove uniqueness, we need to establish that
  \begin{equation}\label{iint}
    \int_{G}\iota\bigl(f(s)\bigr) \,ds=f,
  \end{equation}
  where the equality in \eqref{iint} is meant in
  $M\bigl(\cs(G,\AA)\bigr)$.  Therefore, it suffices to see that
  \begin{equation}\label{eq:2}
    \Bigl(\int\iota\bigl(f(s)\bigr)\,ds\Bigr)
    g=\int_{G}\iota\bigl(f(s)\bigr) g\,ds=f*g
    \quad\text{for all $g\in \sa_c(G,\AA)$.}
  \end{equation}
  Thus we need to establish that the $\cs(G,\AA)$-valued integral in
  the middle of \eqref{eq:2} takes values in (the image of)
  $\sa_{c}(G,\AA)$ in $\cs(G,\AA)$ and coincides with $f*g$.  This can
  be verified almost exactly as in the proof of
  \cite[Lemma~1.108]{danacrossed}.

  Now, if $\rho: C^*(G,\AA)\to M(B)$ is a homomorphism such that
  $\rho\circ\iota = \pi_0$, then by \eqref{iint}, for each
  $f\in\sa_c(G,\AA)$ we must have
  \begin{align*}
    \rho(f) &= \rho\Bigl(\int_G \iota(f(s))\,ds \Bigr) = \int_G
    \rho(\iota(f(s)))\,ds = \int_G \pi_0(f(s))\,ds = \pi(f).
  \end{align*}

  For the converse, let $\pi: C^*(G,\AA)\to M(B)$ be a nondegenerate
  $*$-ho\-mo\-mor\-phism.  By nondegeneracy, $\pi$ extends to a
  strictly continuous homomorphism of $M(C^*(G,\AA))$, so that
  $\pi\circ\iota$ is a strictly continous nondegenerate
  $*$-homomorphism of $\AA$ whose integrated form, by uniqueness, is
  $\pi$.
\end{proof}

\begin{remark}[Modular Differences]
  \label{rem-modular-differences}
  If $p:\AA\to G$ is a Fell bundle over a locally compact
  \emph{group}, then we could just as well have formed the
  $C^{*}$-algebra $\csmw$ by treating $G$ as a groupoid.  (That is, by
  leaving the modular function off the involution.)  To see that
  $\csmw$ and $\csfd$ are naturally isomorphic, we first observe that
  Lemma~\ref{key} and Lemma~\ref{bundle map} remain valid for $\csmw$
  using virtually the same proofs; the only difference is that
  Equations~(\ref{iota}) and (\ref{int}) must be modified to deal with
  the lack of modular function in the involution:
  \begin{gather}
    \label{eq:7}
    \bigl(\iota'(a_{s})f\bigr)(t)=\Delta(s)^{\half}a_{s}f(s^{-1}t)
    \tag*{$(\ref{iota})'$} \quad\text{and} \\
    \pi(f)=\int_{G}\pi_{0}'\bigl(f(s)\bigr) \Delta(s)^{-\half}
    \,ds.\tag*{$(\ref{int})'$}\label{eq:10}
  \end{gather}
  Then notice that there is a $*$-isomorphism $\phi:\samw\to\safd$
  given by $\phi(f)(s)=\Delta(s)^{-\half}f(s)$.  We just need to see
  that $\phi$ is isometric with respect to the universal norm
  $\|\cdot\|_{\text{Gr}}$ on $\csmw$ and $\|\cdot\|$ on $\csfd$.  To
  verify this, let $M$ be a faithful representation of $\csfd$.  Then
  $M$ is the integrated form of $M_{0}:\AA\to B(\HH)$.  But if $L$ is
  the representation of $\csmw$ which is the integrated form of
  $M_{0}$, then
  \begin{equation*}
    \|\phi(f)\|=\|M(\phi(f))\|=\|L(f)\|\le\|f\|_{\text{Gr}}.
  \end{equation*}
  On the other hand, if $L$ is a faithful representation of $\csmw$
  which is the integrated form of $L_{0}$, then we can let $M$ be the
  representation of $\csfd$ that is integrated up from $L_{0}$.  Then
  \begin{equation*}
    \|\phi(f)\|\ge \|M(\phi(f))\|=\|L(f)\|=\|f\|_{\text{Gr}}.
  \end{equation*}
  Thus $\phi$ is isometric.
\end{remark}

\begin{remark}
  \label{rem-more-modular}
  The same comments about modular functions apply to the standard
  group $C^{*}$-algebra and crossed product constructions; that is,
  one can omit the modular function in the definition of the
  involution and arrive at isomorphic algebras. However, you have pay
  for the luxury of modular-free involutions by adding the modular
  function to the integrated form of any representation as in
  \ref{eq:10}.
\end{remark}

\begin{prop}
  \label{prop-Fell-multipliers}
  Let $p:\BB\to \GG$ be a separable Fell bundle over a locally compact
  groupoid $\GG$, and let $\X_{0}$ be a dense subspace of a right
  Hilbert $A$-module $\X$.  Suppose that $L$ is a algebra homomorphism
  of $\sa_{c}(\GG,\BB)$ into the linear operators,
  $\operatorname{Lin}(\X_{0})$, on $\X_{0}$ such that for all
  $x,y\in\X_{0}$
  \begin{enumerate}
  \item $\brip A< L(f)x,y>=\brip A< x,L(f^{*})y>$,
  \item $f\mapsto \brip A< L(f) x,y>$ is continuous in the inductive
    limit topology, and
  \item $\spn\set{L(f)x:\text{$f\in\sa_{c}(\GG,\BB)$ and
        $x\in\X_{0}$}}$ is dense in $\X$.
  \end{enumerate}
  Then $L$ is bounded with respect to the universal $C^*$-norm on
  $\sa_{c}(\GG,\BB)$ and extends to a nondegenerate homomorphism
  $L:\cs(\GG,\BB)\to \mathcal{L}(\X)$.
\end{prop}
\begin{proof}
  This proposition is a consequence of the disintegration result
  \cite[The\-or\-em~4.13]{mw:fell} for Fell bundles.  To see this, let
  $\rho$ be a state on $A$.  Then
  \begin{equation*}
    \ipr(x|y):= \rho\bigl(\rip A<y,x>\bigr)
  \end{equation*}
  is a pre-inner product on $\X_{0}$.  After modding out by the
  subspace $\mathcal{N}$ of vectors of length zero, we get a
  pre-Hilbert space $\HH_{0}:=\X_{0}/\mathcal{N}$ which we view as a
  subspace of its completion $\HH$.  Since
  \begin{equation*}
    \ipr(L(f)x|{L(f)x})=\ipr(x|{L(f^{*}*f)x}),
  \end{equation*}
  it follows from the Cauchy Schwartz inequality that $L(f)$ maps
  $\mathcal{N}$ to itself.  Therefore $L(f)$ defines a linear operator
  $L^{\rho}(f)$ on $\HH_{0}$ via
  $L^{\rho}(f)(x+\mathcal{N})=L(f)x+\mathcal{N}$.  It is clear that
  $L^{\rho}$ defines a pre-representation of $\BB$ on $\HH_{0}$ as in
  \cite[Definition~4.1]{mw:fell}.  Then \cite[Theorem~4.13]{mw:fell}
  implies that
  \begin{equation*}
    \ipr(L(f)x|{L(f)x})\le \|f\|^{2}\ipr(x|x).
  \end{equation*}
  Since this holds for all states $\rho$, we have $\|L(f)\|\le\|f\|$.
  The rest is straightforward.
\end{proof}

\begin{prop}\label{extends}
  Let $\AA$ be a separable Fell bundle over a groupoid $\GG$.  Every
  $*$-homomorphism from $\sa_{c}(\GG,\AA)$ into a $C^*$-algebra
  which
  is continuous from the inductive limit topology into the norm
  topology is bounded for the universal norm, and hence has a unique
  extension to $\cs(\GG,\AA)$.
\end{prop}

\begin{proof}
 Suppose that $\pi:\sa_{c}(\GG,\AA)\to B$ is such a homomorphism,
 and that $\rho: B\to B(\HH)$ is a
 faithful representation of $B$ on a Hilbert space $\HH$. Let
  \[
  \HH_1 = \clspn\set{ \rho\circ\pi(f)\xi : f\in \sa_c(\GG,\AA),
    \xi\in\HH }.
  \]
  Then $f\mapsto \rho\circ\pi(f)|_{\HH_1}$ is a representation of
  $\AA$ on $\HH_1$ in the sense of \cite[Definition~4.7]{mw:fell},
  since the operator norm topology is stronger than the weak operator
  topology.  By \cite[Remark~4.14]{mw:fell},
  \begin{equation*}
    \|\pi(f)\| = \| \rho\circ\pi(f) \| = \| \rho\circ\pi(f)|_{\HH_1} \|
    \leq \| f \|\quad\text{for all $f\in \sa_{c}(\GG,\AA)$.}\qed
  \end{equation*}
  \renewcommand\qed{\relax}
\end{proof}

\section{Product bundles}
\label{product bundles}

If $p:\AA\to G$ is a Fell bundle over a locally compact group $G$,
then the Cartesian product, $\AA\times G$, carries a natural Fell
bundle structure over $G\times G$.  The bundle projection $q:\AA\times G\to
G\times G$ is given by $q(a,t)=(p(a),t)$ and the multiplication and
involution are given by
\[
(a_s,t)(b_r,u)=(a_sb_r,tu) \quad\text{and}\quad
(a_s,t)^*=(a_s^*,t\inv).
\]
(Indeed, the map $(a,t)\mapsto (a,(p(a),t))$ is a bijection of $\AA\times
G$ onto the pull-back Fell bundle $\phi^{*}\AA$ 
--- see Lemma~\ref{lem-pull-back} --- where
$\phi:G\times G\to G$ is the projection onto the first
factor.)

Every section $h\in \sa_c(G\times G,\AA\times G)$ is of the form
\[
h(s,t)=(h_1(s,t),t),
\]
where $h_1\in C_c(G\times G,\AA)$ satisfies $h_1(s,t)\in A_s$
for $s,t\in G$.  For
$f\in \sa_c(G,\AA)$ and $g\in C_c(G)$ we let $f\boxtimes g$ denote the
element of $\sa_c(G\times G,\AA\times G)$ defined by
\[
(f\boxtimes g)(s,t)=(f(s)g(t),t).
\]

\begin{lem}
  \label{basic dense}
  With the above notation, $\spn\set{f\boxtimes g:\text{$f\in
      \sa_c(G,\AA)$ and $g\in C_c(G)$}}$ is inductive-limit dense in
  $\sa_c(G\times G,\AA\times G)$.
\end{lem}

\begin{proof}
  Put $ \SS=\set{f\boxtimes g:f\in \sa_c(G,\AA),g\in C_c(G)}$.  Then
  for each $(s,t)\in G\times G$, $\set{h(s,t):h\in\SS}$ is easily seen
  to be dense in $A_s\times \{t\}$, which is the fibre of the bundle
  $\AA\times G$ over $(s,t)$.  Furthermore if $u,v\in C_{c}(G)$ and
  $u\tensor v$ is the function in $C_{c}(G\times G)$ given by
  $u\tensor v(s,t)=u(s)v(t)$, then $(u\tensor v)h\in\SS$ for all
  $u,v\in C_{c}(G)$ and $h\in\SS$.  Then, because the $u\tensor v$'s
  span an inductive-limit dense subspace of $C_{c}(G\times G)$, a
  straightforward partition of unity argument implies that $\spn\SS$
  is dense as required (see \cite[Proposition~II.14.6 and its
  remark]{fd1} or \cite[Proposition~C.24]{danacrossed}).
\end{proof}

For the study of the coaction associated to a Fell bundle over a group
(specifically, in Section~\ref{coaction crossed product}) we will need
the following slight variation on Lemma~\ref{basic dense}:

\begin{lem}\label{dense}
  Let $\AA\to G$ be a Fell bundle. For $f\in \sa_c(G,\AA)$ and $g\in
  C_c(G)$ define $f\bullet g\in C_c(G\times G,\AA)$ by
  \[
  f\bullet g(s,t)=f(s)g(s\inv t),
  \]
  and define
  \[
  f\star g(s,t)=(f\bullet g(s,t),t).
  \]
  Then $f\star g\in \sa_c(G\times G,\AA\times G)$, and such sections
  have inductive-limit-dense span.
\end{lem}

\begin{proof}
  It is obvious that $f\star g\in \sa_c(G\times G,\AA\times G)$. For
  the second statement, let
  \[
  \SS=\spn\{\,f\star g:\text{$f\in \sa_c(G,\AA)$ and $g\in
    C_c(G)$}\,\}.
  \]
  To show that $\SS$ is dense, we want to invoke a partition of unity
  argument exactly as in Lemma~\ref{basic dense}; thus it suffices to
  establish the following two assertions:
  \begin{enumerate}
  \item For each $(s,t)\in G\times G$, the set $\set{h(s,t):h\in \SS}$
    is dense in $A_s\times \{t\}$;

  \item For each $\kappa,\eta\in C_c(G)$ and $h\in \SS$ we have
    $(\kappa\bullet \eta)h\in \SS$, where similarly to the above we
    define $ \kappa\bullet \eta(s,t)=\kappa(s)\eta(s\inv t)$.

  \end{enumerate}
  (Note that (ii) suffices since the set of functions of the form
  $\kappa\bullet \eta$ have dense span in $C_c(G\times G)$ for the
inductive limit topology,
because this set is the image of the set
  $\set{u\otimes v:u,v\in C_c(G)}$ under the linear homeomorphism
  $\Psi:C_c(G\times G)\to C_c(G\times G)$ defined by
  \[
  \Psi(\p)(s,t)=\p(s,s\inv t),
  \]
  and the functions $u\otimes v$ have dense span in the inductive
  limit topology.)

  For (i), if $a_s\in A_s$ we can choose $f\in\sa_c(G,\AA)$ and $g\in
  C_c(G)$ such that $f(s)=a_s$ and $g(s\inv t)=1$, and then
  \[
  f\star g(s,t)=(a_s,t).
  \]
  For (ii), just observe that
  \[
  (\kappa\bullet \eta)(f\star g)=(\kappa f)\star (\eta g).\qed
  \]
  \renewcommand\qed{\relax}
\end{proof}

\section{Coactions from Fell bundles}
\label{exist}

As mentioned in the introduction, if
$\alpha$ is an action of a locally compact group~$G$ on a $C^*$-algebra~$B$,
then $\AA=B\times G$
has a natural Fell-bundle structure such that
$\cs(G,\AA)\cong B\rtimes_{\alpha}G$.  Then the dual coaction on
$B\rtimes_{\alpha}G$
gives us a coaction on $\cs(G,\AA)$.  In this section, we show that if
$p:\AA\to G$ is any Fell bundle, then $\cs(G,\AA)$ admits a natural
coaction $\delta$ generalizing the dual coaction construction just
described.

\begin{prop}\label{exists}
  Let $\AA$ be a separable Fell bundle over a group $G$.  There is a
  unique coaction $\d$ of $G$ on $C^*(G,\AA)$ such that
  \begin{equation}
    \label{coaction}
    \d(\iota(a_s))=\iota(a_s)\otimes s\quad\text{for $ a_s\in A_s$ and $
      s\in G$.} 
  \end{equation}
\end{prop}

\begin{proof}
  For the proof we will make explicit the canonical map $u: G\to
  M(C^*(G))$.  Consider the map $\d_0: \AA\to M(C^*(G,\AA)\otimes
  C^*(G))$ defined by $\d_0(a_s) = \iota(a_s)\otimes u(s)$.  This
  clearly gives a $*$-homomorphism of $\AA$, and nondegeneracy of
  $\d_0$ follows directly from nondegeneracy of $\iota$.  That $\d_0$
  is strictly continuous follows from strict continuity of $\iota:
  \AA\to M(C^*(G,\AA))$ and $u: G\to M(C^*(G))$.  To see this, let
  $a_{s_i}\to a_s$ in $\AA$, and let $x\in \cs(G,\AA)\otimes
  C^*(G)$. Since $\cs(G,\AA)$ embeds nondegenerately in
  $M(\cs(G,\AA)\otimes C^*(G))$ via $b\mapsto b\otimes 1$, by the
  Hewitt-Cohen factorization theorem we can write $x=(b\otimes 1)y$
  for some $b\in \cs(G,\AA)$ and $y\in \cs(G,\AA)\otimes
  C^*(G)$. Since $\iota(a_{s_i})b\to \iota(a_s)b$ in norm, we have
  $\iota(a_{s_i})b\otimes 1\to \iota(a_s)b\otimes 1$ in norm in
  $M(\cs(G,\AA)\otimes C^*(G))$. Since the map $u:G\to M(C^*(G))$ is
  strictly continuous, and since $a_{s_i}\to a_s$ implies $s_i\to s$
  in $G$, we have $(1\otimes u(s_i))y\to (1\otimes u(s))y$ in norm in
  $\cs(G,\AA)\otimes C^*(G)$. Since multiplication is norm continuous,
  \begin{align*}
    \bigl(\iota(a_{s_i})\otimes u(s_i)\bigr)x
    &=\bigl(\iota(a_{s_i})\otimes u(s_i)\bigr)(b\otimes 1)y
    =\bigl(\iota(a_{s_i})b\otimes 1\bigr) \bigl(1\otimes u(s_i)\bigr)y
    \\
    \intertext{converges in norm to} & \bigl(\iota(a_s)b\otimes
    1\bigr) \bigl(1\otimes u(s)\bigr)y =\bigl(\iota(a_s)\otimes
    u(s)\bigr)x.
  \end{align*}

  Thus Lemma~\ref{bundle map} gives a unique
 nondegenerate
  $*$-ho\-mo\-mor\-phism $\d: C^*(G,\AA)\to M(C^*(G,\AA)\otimes
  C^*(G))$ such that $\d\circ\iota = \d_0$, and by \eqref{int} we have
  \[
  \d(f) = \int_G \iota(f(s))\otimes u(s) \,ds \quad\text{for $f\in
    \sa_c(G,\AA)$.}
  \]

  To see that $\d$ is injective, let $1_G:G\to \C$ be the constant
  function with value $1$, and regard $1_G$ as an element of the
  Fourier-Stieltjes algebra $B(G)=C^*(G)^*$. Then for $f\in
  \sa_c(G,\AA)$ equation~\eqref{iint} and strict continuity of the
  slice map give
  \[
  (\id\otimes 1_G)(\d(f)) =\int (\id\otimes 1_G)(\iota(f(s))\otimes
  u(s))\,ds =\int \iota(f(s))\,ds =f.
  \]
  Thus $(\id\otimes 1_G)\circ \d=\id_{\csga}$ by continuity and
  density, so $\d$ is injective.

  Now if $a_s\in A_s$, then
  \begin{gather*}
    (\d\otimes \id)\circ\d_0(a_s)) =(\d\otimes \id)(\iota(a_s)\otimes
    u(s))
    =\iota(a_s)\otimes u(s)\otimes u(s)\\
    =(\id\otimes \d_G)(\iota(a_s)\otimes u(s)) =(\id\otimes \d_G)\circ
    \d_0(a_s).
  \end{gather*}
  Thus the coaction identity~\eqref{comodule} follows from uniqueness
  in Lemma~\ref{bundle map} together with the usual manipulations with
  vector valued integrals as justified, for example, in
  \cite[Lemma~C.11]{tfb}. 

  Finally, for the nondegeneracy condition \eqref{nondegenerate}, we
  elaborate on the argument sketched in the paragraph preceding
  \cite[Lemma~1.3]{exelng}.  Consider the map $\z_0: \AA\times G\to
  M(\csga\otimes C^*(G))$ defined by $\z_0(a_s,t) = \iota(a_s)\otimes
  u(t)$, where $\AA\times G$ is the Fell bundle over $G\times G$
  defined in Section~\ref{product bundles}.  Arguing as for $\d_0$
  shows that $\z_0$ is a strictly continuous nondegenerate
  $*$-homomorphism, and so Lemma~\ref{bundle map} gives a
  nondegenerate $*$-homomorphism $\z: C^*(G\times G,\AA\times G)\to
  M(C^*(G,\AA)\otimes C^*(G))$ such that $\z\circ\iota = \z_0$.

  In particular, using \eqref{int} and~\eqref{iint} we have, for $f\in
  \sa_c(G,\AA)$ and $g\in C_c(G)$,
  \begin{align*}
    \z(f\boxtimes g)
    &= \int_{G\times G} \z_0\bigl((f\boxtimes g)(s,t)\bigr)\,d(s,t)\\
    &= \int_G \int_G \z_0(f(s)g(t),t) \,ds\,dt\\
    &= \int_G \iota(f(s))\,ds \otimes \int_G g(t)u(t)\,dt\\
    &= f\otimes g,
  \end{align*}
  which implies that $\z$ maps $C^*(G\times G,\AA\times G)$ onto (and
  into) $C^*(G,\AA)\otimes C^*(G)$.

  Similarly, if $f\star g$ is the element of $\sa_c(G\times
  G,\AA\times G)$ defined in Lemma~\ref{dense}, then for $a\otimes
  b\in C^*(G,\AA)\otimes C^*(G)$ we have
  \begin{align*}
    \z(f\star g)(a\otimes b)
    &= \int_{G\times G} \zeta_0(f\star g)(s,t)(a\otimes b)\,d(s,t)\\
    &= \int_G\int_G \zeta_0(f(s)g(s\inv t),t)(a\otimes b)\,dt\,ds\\
    &= \int_G\int_G \iota(f(s))a \otimes g(s\inv t)u(t)b \,dt\,ds\\
    \intertext{which, after $t\mapsto st$, is}
    &= \int_G\int_G \iota(f(s))a \otimes g(t)u(st)b\,dt\,ds\\
    &= \int_G\int_G \bigl(\iota(f(s))\otimes u(s)\bigr)
    \bigl(a\otimes g(t)u(t)b\bigr)\,dt\,ds\\
    &= \Bigl(\int_G \d_0(f(s))\,ds\Bigr)
    \Bigl(\int_G 1\otimes g(t)u(t)\,dt\Bigr)(a\otimes b)\\
    &=\d(f)(1\otimes g)(a\otimes b).
  \end{align*}
  Thus, the multiplier $\d(f)(1\otimes g)$ of $C^*(G,\AA)\otimes
  C^*(G)$ coincides with the image $\zeta(f\star g)$, and the set of
  sections of the form $f\star g$ was shown in Lemma~\ref{dense} to
  have dense span in $\sa_c(G\times G,\AA\times G)$, so the images
  $\z(f\star g)$ have dense span in $C^*(G,\AA)\otimes C^*(G)$.  It
  follows that $\delta$ satisfies the nondegeneracy
  condition~\eqref{nondegenerate}.
\end{proof}

\begin{rem}
  It is clear from the above proof that saturation of the Fell bundle
  $\AA\to G$ is not necessary for Proposition~\ref{exists}.
\end{rem}

\begin{rem}
  Not every coaction is isomorphic to one constructed from a Fell
  bundle as in Proposition~\ref{exists} \cite[Example~2.3(6)]{lprs}.
  For \emph{abelian} $G$, in \cite[Theorem~11.14]{exelspectral} Exel
  effectively characterizes which coactions do arise from Fell bundles
  (modulo the correspondence between coactions of $G$ and actions of
  the Pontryagin dual group $\widehat G$).
\end{rem}

\begin{prop}
  \label{bundle covariant}
  Let $\AA$ be a separable Fell bundle over a group $G$, and let $\d$
  be the coaction of $G$ on $\cs(G,\AA)$ described in
  Proposition~\ref{exists}.  Further let $\pi_0: \AA\to M(B)$ be a
  strictly continuous nondegenerate $*$-homomorphism, with integrated
  form $\pi: C^*(G,\AA)\to M(B)$, and let $\mu: C_0(G)\to M(B)$ be a
  nondegenerate homomorphism.  Then the pair $(\pi,\m)$ is a covariant
  homomorphism of $(\cs(G,\AA),G,\d)$ if and only if
  \begin{equation}
    \label{bundle covariance condition}
    \pi_0(a_s)\m(f)=\m\circ \lt_s(f)\pi_0(a_s)
    \quad\text{for $s\in G$, $a_s\in A_s$ and $f\in C_0(G)$,}
  \end{equation}
  where $\lt$ is the action of $G$ on $C_0(G)$ by left translation:
  $\lt_s(f)(t)=f(s\inv t)$.
\end{prop}

\begin{proof}
  First assume that $(\pi,\m)$ is covariant.  Because $B(G)\cap
  C_0(G)$ is dense in $C_0(G)$, it suffices to verify \eqref{bundle
    covariance condition} for $f\in B(G)$.  So fix $f\in B(G)$, and
  put $g=\lt_s(f)\in B(G)$.  By~\cite[Proposition~A.34]{enchilada}, we
  have
  \[
  (\id_{B}\tensor g)\bigl((\mu\tensor\id)(w_{G})\bigr) =\mu(g)
  \]
  where $\id\otimes g:M(C_0(G)\otimes C^*(G))\to M(C_0(G))$ denotes
  the slice map.  Then
  \begin{align*}
    \m\circ\lt_s(f)\pi_0(a_s) &=\m(g)\pi_0(a_s)\\
    &=(\id_{B}\otimes
    g)\bigl((\m\otimes\id)(w_G)\bigr)\pi_0(a_s) \\
    \intertext{which, by \cite[Lemma~A.30]{enchilada}, is}
    &=(\id_{B}\otimes g)\bigl((\m\otimes\id)(w_G)(\pi(\iota(a_s))
    \otimes 1)\bigr)
    \\
    \intertext{which, by the covariance condition \eqref{covariant},
      is} &=(\id_{B}\otimes
    g)\bigl((\pi\otimes\id)(\d(\iota(a_s)))(\m\otimes\id)(w_G)\bigr)
    \\
    &=(\id_{B}\otimes g)\bigl((\pi(\iota(a_s))\otimes
    u(s))(\m\otimes\id)(w_G)\bigr) \\
    &=(\id_{B}\otimes g)\bigl((\pi_0(a_s)\otimes
    1)(\m\otimes\id)\bigl((1\otimes u(s))w_G\bigr)\bigr)
    \\
    \intertext{which, after applying \cite[Lemma~A.30]{enchilada} and
      writing $(\lt_{s\inv}\otimes\id)(w_G)$ for the multiplier
      $r\mapsto u(sr)$, is} &=\pi_0(a_s)(\id_{B}\otimes
    g)\bigl((\m\otimes\id)\bigl((\lt_{s\inv}\otimes\id)(w_G)\bigr)\bigr)
    \\
    \intertext{which, since $(\mu\tensor\id)\circ (\lt_{s^{-1}}\tensor
      \id) = \mu\circ\lt_{s^{-1}}\tensor \id$ as a nondegenerate
      homomorphism of $C_{0}(G)\tensor C^{*}(G)$ into $M(B)\tensor
      C^{*}(G)\subseteq M(B\tensor C^{*}(G))$, is}
    &=\pi_0(a_s)(\id_{B}\otimes
    g)\bigl((\m\circ\lt_{s\inv}\otimes\id)(w_G)\bigr)
    \\
    \intertext{ which, by \cite[Proposition~A.34]{enchilada}, is }
    &=\pi_0(a_s)\m\circ\lt_{s\inv}(g) \\
    &=\pi_0(a_s)\m(f).
  \end{align*}

  Conversely, the above computation
  can be rearranged to show that, if \eqref{bundle covariance
    condition} holds, then
  \[
  (\id_{B}\otimes g)\bigl((\m\otimes\id)(w_G)(\pi(\iota(a_s)) \otimes
  1)\bigr) = (\id_{B}\otimes
  g)\bigl((\pi\otimes\id)(\d(\iota(a_s)))(\m\otimes\id)(w_G)\bigr)
  \]
  for every $g\in B(G)$.  Since slicing by elements of $B(G)$
  separates points in $M(\csga\otimes C^*(G))$, it follows that the
  covariance condition \eqref{covariant} holds for every $a$ of the
  form $\iota(a_s)$, which then implies (by Lemma~\ref{bundle map})
  that it holds for every element of $C^*(G,\AA)$.
\end{proof}

We include the following proposition since it might be useful
elsewhere, although we will not need it in the present paper.

\begin{prop}
  If $\alpha$ is an action of a group $G$ on a $C^*$-algebra $B$, and
  $\AA\to G$ is the associated semidirect-product Fell bundle,
then the isomorphism
  \[
  B\rtimes_\alpha G \cong C^*(G,\AA)
  \]
  carries the dual coaction $\hat\alpha$ to the coaction $\delta$ of
  $G$ on $C^*(G,\AA)$ described in Proposition~\ref{exists}.
\end{prop}

\begin{proof}
  We recall that the isomorphism $\theta: B\rtimes_\alpha G\to
  C^*(G,\AA)$ is characterized on generators by
  \[
  \theta(i_B(b)i_G(f)) = \int_G f(s) \iota(b,s)\,ds \quad\text{for
    $b\in B$ and $ f\in C_c(G)$}
  \]
  (which follows from \cite[\S VIII.5.7]{fd2}).  Thus,
  \begin{align*}
    \delta\circ\theta(i_B(b)i_G(f))
    &= \int_G f(s) \delta(\iota(b,s))\,ds\\
    &= \int_G f(s) \iota(b,s)\otimes s\,ds\\
    &= \int_G f(s) (\theta\otimes\id)(i_B(b)i_G(s)\otimes s)\,ds\\
    &= \int_G f(s) (\theta\otimes\id)\circ\hat\alpha(i_B(b)i_G(s))\,ds\\
    &= (\theta\otimes\id)\circ\hat\alpha(i_B(b)i_G(f)).\qed
  \end{align*}
  \renewcommand\qed{\relax}
\end{proof}

\section{Transformation bundles}
\label{transformation bundles}

Having defined a coaction $\delta$ on the $C^{*}$-algebra $\cs(G,\AA)$
of a Fell bundle over a group, an obvious next step is to consider the
corresponding crossed product.  In the next section, we will show that
$\cs(G,\AA)\rtimes_{\delta}G$ is isomorphic to the $C^{*}$-algebra of
a Fell bundle over a groupoid.  The purpose of this short section is
to describe that groupoid and Fell bundle.

Let $G$ be a locally compact group, and let $G\times_{\lt} G$ denote
the transformation groupoid associated to the action $\lt$ of $G$ on
itself by left translation, with multiplication and inverse
\[
(s,tr)(t,r)=(st,r) \quad\text{and}\quad (s,t)\inv = (s\inv,st) \quad
\text{for $s,t,r\in G$.}
\]
Note that the unit space is $(G\times_{\lt} G)^0=\{e\}\times G$, and
the range and source maps are given by
\[
r(s,t)=(e,st)\quad\text{and}\quad s(s,t)=(e,t).
\]
It it not hard to check that we get a left Haar system on
$G\times_{\lt}G$ via
\[
\int_{G\times_\lt G} f(u,v)\,d\l^{r(s,t)}(u,v) = \int_G f(u,u\inv
st)\,du \quad\text{for $f\in C_c(G\times_\lt G)$.}
\]

Now let $\AA\to G$ be a Fell bundle over the locally compact group
$G$.  The map $\phi:(s,t)\mapsto s$ is a groupoid homomorphism of
$G\times_{\lt}G$ onto the group $G$.  The pull-back Fell bundle $\phi^{*}\AA$ 
(see Lemma~\ref{lem-pull-back}) will be called the
\emph{transformation Fell bundle} $\AA\times_{\lt} G\to G\times_{\lt}
G$.  We will use the bijection $(a_s,(s,t))\mapsto (a_s,t)$ to identify
the total space of $\AA\times_{\lt}G$ with the Cartesian product
$\AA\times G$. Then the multiplication is
\[
(a_s,tr)(b_t,r)=(a_sb_t,r)\quad\text{for $s,t,r\in G$, $a_s\in A_s$
  and $b_t\in A_t$,}
\]
and the involution is
\[
(a_s,t)^*=(a_s^*,st).
\]
For future reference, the convolution in
$\sa_c(G\times_{\lt}G,\AA\times_{\lt}G)$ is given by
\begin{align}
  (h*k)(s,t) 
  &=\int_G h(u,u\inv st)k(u\inv s,t)\,du\label{eq:3}
\end{align}
and the involution by
\begin{equation}
  h^*(s,t) =h\bigl((s,t)\inv\bigr)^* =h(s\inv,st)^*.\label{eq:4}
\end{equation}

Note that every $h\in\sa_{c}(\gltg,\altg)$ is of the form
\begin{equation*}
  h(s,t)=\bigl(h_{1}(s,t),t\bigr)
\end{equation*}
for a continuous function $h_{1}:\gltg\to\AA$ with $h_1(s,t)\in A_{s}$.

\section{Coaction crossed product}
\label{coaction crossed product}

Our purpose in this section is to prove the following:

\begin{thm}
  \label{coaction crossed product isomorphism}
  Let $\AA$ be a separable Fell bundle over a group $G$, and let $\d$
  be the associated coaction on $\cs(G,\AA)$ described in
  Proposition~\ref{exists}.  If $q:\altg\to\gltg$ is the
  transformation Fell bundle constructed in the preceding section,
  then there is an isomorphism
  \[
  \t:\cs(G,\AA)\rtimes_\d G\to \cs(G\times_{\lt}G,\AA\times_{\lt} G)
  \]
  such that
  \begin{equation}\label{theta}
    \t\bigl(j_{\cs(G,\AA)}(f)j_{G}(g)\bigr)=(\D^{\half}f)\boxtimes g
    \quad\text{for $f\in \sa_c(G,\AA)$ and $ g\in C_c(G)$,}
  \end{equation}
  where $(\D^{\half}f)\boxtimes g\in \sa_c(G\times_{\lt}G,
  \AA\times_{\lt} G)$ is defined by $((\D^{\half}f)\boxtimes g)(s,t) =
  (\D(s)^{\half}f(s)g(t),t)$.
\end{thm}
\begin{remark}
  For $G$ discrete, this is a special case of
  \cite[Corollary~2.8]{eq:full}.
\end{remark}

\begin{proof}
  We will obtain $\theta$ as the integrated form of a covariant
  homomorphism $(\t_\AA,\t_G)$ of $(\cs(G,\AA),G,\d)$ into
  $M(\cs(G\times_{\lt}G, \AA\times_{\lt} G))$ such that
  \begin{equation}
    \label{theta2}
    \t_\AA(f)\t_G(g)=(\D^{\half}f)\boxtimes g\in \sa_c(G\times_{\lt}G,
    \AA\times_{\lt} G) 
  \end{equation}
  for $f\in \sa_c(G,\AA)$ and $g\in C_c(G)$.  It will follow that
  $\t=\t_\AA\rtimes\t_G$ maps $A\rtimes_\d G$ into $C^*(G\times_\lt
  G,\AA\times_{\lt} G)$, satisfies~\eqref{theta}, and is surjective
  because $\set{f\boxtimes g: f\in \sa_c(G,\AA),g\in C_c(G)}$ has
  inductive-limit-dense span in $\sa_c(G\times_\lt G,\AA\times_{\lt}
  G)$.  We will show that $\theta$ is injective by finding a
  representation $\Pi$ of $C^*(G\times_\lt G, \AA\times_\lt G)$ such
  that $\Pi\circ\theta$ is a faithful regular representation of
  $C^*(G,\AA)\rtimes_\d G$.

  We will obtain $\t_{\AA}:\cs(G,\AA)\to
  M\bigl(\cs(G\times_{\lt}G,\AA\times_{\lt}G)\bigr)$ as the integrated
  form of a $*$-homomorphism $\tz:\AA\to
  M\bigl(\cs(G\times_{\lt}G,\AA\times_{\lt}G)\bigr)$ (as in
  Lemma~\ref{bundle map}).
Given $a_{s}\in A_{s}$, we define an
  operator $\tz(a_s)$ on $\sa_{c}(G\times_{\lt}G,\AA\times_{\lt}G)$ by%
\footnote{The operator $\tz(a_{s})$ defined in
    \eqref{eq:5} is analogous to $\iota(a_{s})$ defined in
    Lemma~\ref{key}.  The modular function appearing in its definition
    is required to make $\tz$ $*$-preserving. It is necessary here
    because there is no modular function in the involution in
    $\sa_{c}(\gltg,\altg)$.}
  \begin{equation}
    \label{eq:5}
    (\tz(a_s)h)(t,r) =\bigl(a_sh_1(s\inv t,r)\D(s)^{\half},r\bigr).
  \end{equation}
Then it is straightforward to verify that
$\tz(a_s)\tz(a_t) = \tz(a_s a_t)$.
Moreover,
if $h,k\in \sa_c(G\times_\lt G,\AA\times_\lt G)$, we have
  \begin{align*}
    \brip< &\tz(a_s)h,k> (t,r)
    = \bigl( (\tz(a_s)h)^**k \bigr)(t,r)\\
    \intertext{which, in view of the formula for convolution given by
      \eqref{eq:3}, is}
    &= \int_G (\tz(a_s)h)^*(u,u\inv tr)\, k(u\inv t,r)\,du\\
    \intertext{which, using the formula for the involution given by
      \eqref{eq:4}, is}
    &= \int_G (\tz(a_s)h)(u\inv,tr)^*\, k(u\inv t,r)\, du\\
    &= \int_G \bigl( a_s h_1(s\inv u\inv,tr)\D(s)^{\half},tr \bigr)^*
    \, \bigl( k_1(u\inv t,r),r \bigr)\, du\\
    &= \int_G \bigl( h_1(s\inv u\inv,tr)^* a_s^*\D(s)^{\half}, u\inv
    tr \bigr)\,
    \bigl( k_1(u\inv t,r),r \bigr)\, du\\
    &= \int_G \bigl( h_1(s\inv u\inv,tr)^* a_s^* k_1(u\inv t,r),
    r\bigr)\, \D(s)^{\half}du\\
    \intertext{which, after sending $u\mapsto us^{-1}$, is} &= \int_G
    \bigl( h_1(u\inv,tr)^* a_s^* k_1(su\inv t,r),r
    \bigr)\,\D(s)^{-\half}du\\
    &= \int_G \bigl( h_1(u\inv,tr)^*, u\inv tr\bigr)
    \bigl( a_s^*k_1(su\inv t,r)\D(s)^{-\half},r\bigr)\,du\\
    &= \int_G \bigl(h_1(u\inv,tr),tr\bigr)^*
    \bigl(a_s^*k_1(su\inv t,r)\D(s)^{-\half},r\bigr)\,du\\
    &= \int_G h(u\inv,tr)^* (\tz(a_s^*)k)(u\inv t,r)\, du\\
    &= \int_G h^*(u,u\inv tr) (\tz(a_s^*)k)(u\inv t,r)\, du\\
    &= \bigl( h^**(\tz(a_s^*)k)\bigr)(t,r)\\
    &= \brip< h,\tz(a_s^*)k >(t,r).
  \end{align*}
If we choose $b_{e}\in\tA_{e}$ such that
$\|a_{s}\|^{2}1_{A_{e}}-a_{s}^{*}a_{s}=b_{e}^{*}b_{e}$, 
then since~$\tz$ makes sense and is multiplicative
on $\tA_{e}$, and since
the preceding computation certainly holds for
  $b_{e}\in\tA_{e}$, we see that
  \begin{align*}
    \|a_{s}\|^{2}\brip< h,h> -\brip< \tz(a_{s})h,\tz(a_{s})h> &=
    \brip<
    \tz(\|a_{s}\|^{2}1_{A_{e}} - a_{s}^{*}a_{s})h,h> \\
    &=\brip< \tz(b_{e})h,\tz(b_{e})h> \ge 0
  \end{align*}
for all $h\in \sa_{c}(G\times_{\lt}G,\AA\times_{\lt}G)$. 
Thus $\tz(a_{s})$ extends to a bounded adjointable operator on
  $\cs(G\times_{\lt}G,\AA\times_{\lt}G)$ and we get a
$*$-homomorphism
$\tz:\AA\to M(C^*(G\times_\lt G,\AA\times_{\lt} G))$.  

We need to show that
$\tz$ is strictly continuous and nondegenerate.
For nondegeneracy, let $\{e_i\}$ be an approximate identity in
  $A_e$. It suffices to show that if $h\in \sa_c(G\times_{\lt}G,
  \AA\times_{\lt} G)$ then $\tz(e_i)h\to h$ in the inductive limit
  topology.\footnote{This could be proved using
    \cite[Lemma~8.1]{mw:fell}.  However, the proof of that lemma given
    in \cite{mw:fell} is incorrect.  Fortunately, it can be fixed
    along the same lines as presented here.}  Notice that
  \[
  \tz(e_i)h(r,t)=(e_ih_1(r,t),t).
  \]
  Since each $A_{r}$ is an $A_{e}$\,--\,$A_{e}$-imprimitivity bimodule,
  $e_ih_{1}(r,t)\to h_{1}(r,t)$ for any $(r,t)\in\gltg$.  Fix
  $\epsilon>0$.  Since $a\mapsto \|a\|$ is continuous on $\AA$, we can
  cover $\supp h_{1}$ with open sets $V_{1},\dots,V_{n}$ and find
  $a_{j}\in A_{e}$ such that
  \begin{equation*}
    \|a_{j}h_{1}(r,t)-h_{1}(r,t)\|<\frac\epsilon3\quad\text{for all
      $(r,t)\in V_{j}$.}
  \end{equation*}
  Let $\set{\phi_{j}}\subseteq C_{c}^{+}(\gltg)$ be such that
  $\supp\phi_{j}\subseteq V_{j}$ and 

$\sum_{j}\phi_{j}(r,t)\leq1$ for all
  $(r,t)$, with equality for $(r,t)\in\supp h_{1}$.
 Define $a\in
  C_{c}(\gltg,A_{e})$ by
  \begin{equation*}
    a(r,t)=\sum_{j}\phi_{j}(r,t)a_{j}.
  \end{equation*}
  Then
  \begin{equation*}
    \|a(r,t)h_{1}(r,t)-h_{1}(r,t)\|<\frac\epsilon3\quad\text{for all $(r,t)$.}
  \end{equation*}
  Clearly, there is an $i_{0}$ such that $i\ge i_{0}$ implies that
  \begin{equation*}
    \|e_{i}a(r,t)-a(r,t)\|<\frac\epsilon{3(\|h_{1}\|_{\infty}+1)}\quad\text{for
      all $(r,t)$.}
  \end{equation*}
  Since $\|e_{i}\|\le1$ for all $i$, we see that $i\ge i_{0}$ implies
  \begin{align*}
    \|\tz(e_{i})&h(r,t)-h(r,t)\|= \|e_{i}h_{1}(r,t)-h_{1}(r,t)\| \\
    &\le
    \|e_{i}h_{1}(r,t)-e_{i}a(r,t)h_{1}(r,t)\|+\|e_{1}a(r,t)h_{1}(r,t)
    - a(r,t)h_{1}(r,t)\| \\
    &\hskip1in +
    \|a(r,t)h_{1}(r,t)-h_{1}(r,t)\| \\
    &\le 2 \|h_{1}(r,t)-a(r,t)h_{1}(r,t)\|
    +\|e_{i}a(r,t)-a(r,t)\|\|h_{1}\|_{\infty} \\
    &<\frac\epsilon3+\frac\epsilon3+\frac\epsilon3=\epsilon.
  \end{align*}
  Therefore $\tz(e_i)h\to h$ uniformly,
so since $\supp \tz(e_i)h = \supp h$ for all $i$,
we have $\tz(e_i)h\to h$ in the inductive limit topology,
as desired.

  Finally, for strict continuity we note that our separability
  assumption on $p:\AA\to G$ guarantees that $\AA$ is second countable
  \cite[Proposition~II.13.21]{fd1}.  Thus, it suffices to show that
$\tz$ takes convergent sequences to strictly convergent sequences.

So suppose 
  $\set{a_{i}}$ is a sequence converging to $a$ in $\AA$. 
Let $s=p(a)$, and 
   for each~$i$, let $s_i=p(a_i)$;  so $s_i\to s$ in $G$.
  Since $\set{a_{i}}$ must lie in a norm-bounded subset of $\AA$, the
  image $\set{\tz(a_{i})}$ is a bounded sequence in
  $M(\cs(G\times_{\lt}G, \AA\times_{\lt}G))$.  Thus it suffices to
  show that $\tz(a_{i})\to\tz(a)$ $*$-strongly
  \cite[Proposition~C.7]{tfb}.  Since $a_{i}^{*}\to a^{*}$ and $\tz$
  is $*$-preserving, it suffices to show strong convergence.  Since
  $\set{\tz(a_{i})}$ is bounded, it suffices to show that
  $\tz(a_{i})h\to \tz(a)h$ in the inductive limit topology for each
  $h\in \sa_{c}(G\times_{\lt}G,\AA\times_{\lt}G)$.

  We can replace $\{a_i\}$ by a subsequence (keeping the same
  notation) such that the $s_i$'s lie in a fixed compact neighborhood
  of $s$. Then the supports of the $\tz(a_i)h$'s all lie in a fixed
  compact set, so it suffices to show that $\tz(a_i)h\to \tz(a)h$
  uniformly.
  If not, then there are $(r_{i},t_{i})$, all lying in a compact
  subset of $\gltg$, and an $\epsilon>0$ such that
  \begin{equation}
    \label{eq:11}
    \|\tz(a_{i})h(r_{i},t_{i})-\tz(a)h(r_{i},t_{i})\|\ge\epsilon.
  \end{equation}
  Of course, we can pass to a subsequence, relabel, and assume that
  $(r_{i},t_{i})\to (r,t)$.  But the left-hand side of \eqref{eq:11}
  equals
  \begin{equation}
    \label{eq:12}
    \|a_{i}h_{1}(s_{i}^{-1}r_{i},t_{i})-ah_{1}(s^{-1}r_{i},t_{i})\|.
  \end{equation}
  Since $(a_{i},h_{1}(s_{i}^{-1}r_{i},t_{i}))$ and
  $(a,h_{1}(s^{-1}r_{i},t_{i}))$ both converge to
  $(a,h_{1}(s^{-1}r,t))$ in $\AA\times \AA$, and since multiplication
  is continuous from $\AA\times\AA\to \AA$, it follows that
  $a_{i}h_{1}(s_{i}^{-1}r_{i},t_{i})-ah_{1}(s^{-1}r_{i},t_{i})$ tends
  to $0_{A_{r}}$ in $\AA$.  Therefore, \eqref{eq:12} tends to zero,
  and this contradicts \eqref{eq:11}.  Thus $\tz$ is strictly
  continuous.

  Having dealt with $\tz$, we turn to the definition of $\t_{G}$.  For
  $f\in C_0(G)$ and $h\in \sa_c(G\times_{\lt}G,\AA\times_{\lt} G)$
  define
  \begin{equation}
    \label{eq:6}
    (\t_G(f)h)(s,t)=(f(st)h_1(s,t),t).
  \end{equation}
  We note that \eqref{eq:6} makes perfectly good sense for $f\in
  C_{0}(G)^{\sim}$,
and then $\t_G(fg) = \t_G(f)\t_G(g)$ for $f,g\in C_{0}(G)^{\sim}$.
Another computation shows that
  \begin{equation*}
    \brip< \t_{G}(f)h,k> = \brip< h,\t_{G}(\bar f)k>
  \end{equation*}
for all such $f$.  
Writing
$\|f\|_{\infty}^{2}-\bar f f = \bar g g$ for some $g\in C_{0}(G)^{\sim}$,
we thus have
  \begin{equation*}
    \|f\|_{\infty}^{2}\brip <h,h> - \brip<
    \t_{G}(f)h,\t_{G}(f)h> = \brip<
    \t_{G}(g)h,\t_{G}(g)h> \ge 0
  \end{equation*}
for all $h$.  
  Therefore $\t_{G}(f)$ is bounded and we get a $*$-homomorphism of
  $C_{0}(G)$ into $M\bigl(\cs(G\times_{\lt}G,\AA\rtimes_{\lt
    G})\bigr)$.\footnote{In fact, modulo the obvious identification of
    $G$ with $(G\times_{\lt}G)^{(0)}$, $\t_{G}$ is just the natural
    map of $C_{0}(\GG^{(0)})$ into the multiplier algebra of the
    $C^{*}$-algebra $\cs(\GG,\BB)$ of a Fell bundle over a groupoid
    $\GG$.}

  We let $\t_{\AA}$ be the integrated form of $\tz$ (see
  Lemma~\ref{bundle map}).  To see that $(\t_{\AA},\t_{G})$ is
  covariant, we will use Proposition~\ref{bundle covariant}.  For
  $a_s\in A_s$, $f\in C_c(G)$, $h\in \sa_c(G\times G,\AA\times G)$,
  and $r,t\in G$ we have
  \begin{align*}
    \bigl(\tz(a_s)\t_G(f)h\bigr)(r,t) &=\bigl(a_s(\t_G(f)h)_1(s\inv
    r,t)\D(s)^{\half},r\bigr) \\&=\bigl(a_sf(s\inv rt)h_1(s\inv
    r,t)\D(s)^{\half},r\bigr) \\&=\bigl(\lt_s(f)(rt)a_sh_1(s\inv
    r,t)\D(s)^{\half},r\bigr)
    \\&=\bigl(\lt_s(f)(rt)(\tz(a_s)h)_1(r,t),r\bigr)
    \\&=\bigl(\t_G\circ\lt_s(f)\tz(a_s)h\bigr)(r,t).
  \end{align*}

  To verify \eqref{theta2}, for $h\in \sa_c(G\times_\lt
  G,\AA\times_{\lt}G)$ we have
  \begin{align*}
    \bigl(\t_\AA(f)\t_G(g)h\bigr)(s,t) &=\int
    \bigl(\tz(f(r))\t_G(g)h\bigr)(s,t)\,dr \\&=\int
    \bigl(f(r)(\t_G(g)h)_1(r\inv s,t)\D(r)^{\half},t\bigr)\,dr
    \\&=\int \bigl(f(r)g(r\inv st)h_1(r\inv
    s,t)\D(r)^{\half},t\bigr)\,dr \\&=\int
    \bigl(f(r)\D(r)^{\half}g(r\inv st),r\inv st\bigr)\bigl(h_1(r\inv
    s,t),t\bigr)\,dr \\&=\int \bigl((\D^{\half}f)\boxtimes
    g\bigr)(r,r\inv st)h(r\inv s,t)\,dr
    \\&=\bigl(\bigl((\D^{\half}f)\boxtimes g\bigr)*h\bigr)(s,t).
  \end{align*}

As outlined at the start of the proof, 
it follows from the above that the integrated form
$\t=\t_\AA\rtimes\t_G$ maps $A\rtimes_\d G$ (into and) onto $C^*(G\times_\lt
  G,\AA\times_{\lt} G)$.  To show that $\t$ is faithful,
we will now construct a
  representation $\Pi$ of $C^*(G\times_\lt G,\AA\times_\lt G)$ such
  that $\Pi\circ\theta$ is the regular representation
$\L=(\pi_\AA\otimes\l)\circ \d\rtimes (1\otimes M)$
associated to a faithful representation $\pi_\AA$ of $\cs(G,\AA)$.
This will
  suffice since $\L$ is faithful by \cite[Remark~A.43(3)]{enchilada}.

  So let $\pi_\AA$ be a faithful nondegenerate representation of
  $\cs(G,\AA)$ on a Hilbert space $\HH$.
  Of course, $\pi_{\AA}$ is the integrated form of a representation
  $\piz$ of $\AA$, by Lemma~\ref{bundle map}.  For $h\in
  \sa_c(G\times_{\lt}G, \AA\times_{\lt} G)$ and $\xi\in
  C_c(G,\HH)\subseteq \HH\otimes L^2(G)$, define $\Pi_0(h)\xi: G\to \HH$
  by
  \begin{equation}\label{Pi}
    \bigl(\Pi_0(h)\xi\bigr)(t)=\int_G \piz\bigl(h_1(s,s\inv t)\bigr)
    \xi(s\inv t)\D(s)^{-\half}\,ds;
  \end{equation}
  the integrand is in $C_c(G\times G,\HH)$, so \eqref{Pi} does define
  a vector in $\HH$, and $\Pi_0(h)\xi\in C_c(G,\HH)$.  It follows that
  \eqref{Pi} defines a linear operator $\Pi_0(h)$ on the dense
  subspace $C_c(G,\HH)$ of $\HH\otimes L^2(G)$.

  By \cite[Theorem~4.13]{mw:fell}, to show that $\Pi_0$ extends to a
  representation $\Pi:C^*(G\times_\lt G,\AA\times_{\lt} G)\to
  B(\HH\otimes L^2(G))$, it suffices to show that $\Pi_0$ is a
  \emph{pre-representation} of $\AA\times_{\lt} G$ on
  $C_c(G,\HH)$. Recall from \cite[Definition~4.1]{mw:fell} that to say
  that $\Pi_0$ is a pre-representation means that
  $\Pi_0:\sa_c(G\times_{\lt}G, \AA\times_{\lt} G)\to \lin(C_c(G,\HH))$
  (where $\lin(C_c(G,\HH))$ denotes the algebra of all linear
  operators on the vector space $C_c(G,\HH)$) is an algebra
  homomorphism such that for all $\xi,\eta\in C_c(G,\HH)$:
  \begin{enumerate}
  \item $h\mapsto \brip<\Pi_0(h)\xi,\eta>$ is continuous in the
    inductive limit topology;
  \item $\brip<\Pi_0(h)\xi,\eta>=\brip<\xi,\Pi_0(h^*)\eta>$; and
  \item $\Pi_0(\sa_c(G\times_\lt G,\AA\times_{\lt} G))C_c(G,\HH)$ has
    dense span in $\HH\otimes L^2(G)$.
  \end{enumerate}
  $\Pi_0$ is obviously linear; we verify that it is multiplicative:
  for $f,g\in \sa_c(G\times_{\lt}G,\AA\times_{\lt} G)$ and $\xi\in
  C_c(G,\HH)$ we have
  \begin{align*}
    \bigl(\Pi_0&(f*g)\xi\bigr)(t) =\int_G\piz\bigl((f*g)_1(s,s\inv
    t)\bigr)\xi(s\inv
    t)\D(s)^{-\half}\,ds \\
    &=\int_G\int_G\piz\bigl(f_1(r,r\inv t)g_1(r\inv s,s\inv
    t)\bigr)\xi(s\inv t)\D(s)^{-\half}\,ds\,dr
    \\
    \intertext{which, after $s\mapsto rs$, is}
    &=\int_G\int_G\piz\bigl(f_1(r,r\inv t)g_1(s,s\inv r\inv
    t)\bigr)\xi(s\inv r\inv t)\D(rs)^{-\half}\,ds\,dr
    \\
    &=\int_G\piz\bigl(f_1(r,r\inv t)\bigr) \\
    &\hskip1in \Bigl(\int_G\piz\bigl(g_1(s,s\inv r\inv
    t)\bigr)\xi(s\inv r\inv t) \D(s)^{-\half}\,ds\Bigr)
    \D(r)^{-\half}\,dr
    \\
    &=\bigl(\Pi_0(f)\Pi_0(g)\xi\bigr)(t).
  \end{align*}

  For (i), it suffices to show that if $K\subseteq G\times G$ is compact
  and $\{h_n\}$ is a sequence converging uniformly to $0$ in
  $\G_K(G\times_\lt G,\AA\times_{\lt} G)$ then
  \[
  \brip<\Pi_0(h_n)\xi,\eta>\to 0\quad\text{for all $\xi,\eta\in
    C_c(G,\HH)$.}
  \]
  We have
  \begin{align*}
    \brip<\Pi_0(h_n)\xi,\eta>
    &=\int_G\brip<\bigl(\Pi_0(h_n)\xi\bigr)(t),\eta(t)>\,dt
    \\
    &=\int_G\int_G \brip<\piz({h_n(s,s\inv t))\xi(s\inv
      t)},\eta(t)>\,ds\,dt,
  \end{align*}
  which converges to $0$ since the integrands converge uniformly to
  $0$ and the integration is over a compact set.

  For (ii) we have
  \begin{align*}
    \brip<\Pi_0(h)\xi,\eta>
    &=\int_G\brip<\bigl(\Pi_0(h)\xi\bigr)(t),\eta(t)>\,dt
    \\
    &=\int_G\int_G\brip<{\piz(h_1(s,s\inv t))\xi(s\inv t)},\eta(t)>
    \D(s)^{-\half}\,dt\,ds
    \\
    &=\int_G\int_G\brip<\xi(s\inv t),\piz(h_1(s,s\inv t)^*)\eta(t)>
    \D(s)^{-\half}\,dt\,ds
    \\
    \intertext{which, after $t\mapsto st$, is}
    &=\int_G\int_G\bigl\<\xi(t),\piz(h_1(s,t)^*)\eta(st)\bigr\>
    \D(s)^{-\half}\,dt\,ds
    \\
    \intertext{which, after $s\mapsto s\inv$, is}
    &=\int_G\int_G\brip<\xi(t),\piz(h_1(s\inv,t)^*)\eta(s\inv t)>
    \D(s)^{-\half}\,dt\,ds
    \\
    &=\int_G\int_G\brip<\xi(t),\piz((h^*)_1(s,s\inv t))\eta(s\inv t)>
    \D(s)^{-\half}\,ds\,dt
    \\
    &=\int_G\brip<\xi(t),\bigl(\Pi_0(h^*)\eta\bigr)(t)>\,dt
    \\
    &=\brip<\xi,\Pi_0(h^*)\eta>.
  \end{align*}

  For (iii), it suffices to show that for $f\in \sa_c(G,\AA)$ and
  $g\in C_c(G)$ we have
  \[
  \Pi_0((\D^{\half}f)\boxtimes g) =
  (\pi_\AA\otimes\l)\circ\d(f)(1\otimes M_g),
  \]
  because the ranges of the operators on the right-hand side have
  dense span in $\HH\otimes L^2(G)$ since the regular representation
  of $A\times_\d G$ is nondegenerate. For $\xi\in C_c(G,\HH)$ we have
  \begin{align*}
    \Bigl(\Pi_0\bigl((\D^{\half}f)\boxtimes g\bigr)\xi\Bigr)(t)
    &=\int_G\piz\bigl((f\boxtimes g)_{1}(s,s\inv t)\bigr)\xi(s\inv
    t)\D(s)^{-\half}\,ds \\
    &=\int_G\piz(f(s))g(s\inv t)\xi(s\inv
    t)\,ds \\
    &=\int_G\piz(f(s))(M_g\xi)(s\inv t)\,ds
    \\
    &=\int_G\piz(f(s))(\l_sM_g\xi)(t)\,ds
    \\
    &=\int_G\bigl(\bigl(\piz(f(s))\otimes
    \l_sM_g\bigr)\xi\bigr)(t)\,ds
    \\
    &=\int_G\bigl(\bigl(\pi_\AA(\iota(f(s)))\otimes
    \l_sM_g\bigr)\xi\bigr)(t)\,ds
    \\
    &=\int_G\bigl(\bigl(\pi_\AA\otimes\l\bigr)\bigl(\iota(f(s))\otimes
    u(s)\bigr) \bigl(1\otimes M_g\bigr)\xi\bigr)(t)\,ds
    \\
    &=\int_G\bigl((\pi_\AA\otimes\l)\circ\d(\iota(f(s)))(1\otimes
    M_g)\xi\bigr)(t)\,ds
    \\
    &=\bigl((\pi_\AA\otimes\l)\circ\d(f)(1\otimes M_g)\xi\bigr)(t).
  \end{align*}

  As we explained above, we now can conclude that $\Pi_0$ extends
  uniquely to a nondegenerate representation $\Pi$ of
  $\cs(G\times_{\lt}G,\AA\times_{\lt} G)$, and then the above
  calculation verifies that $\Pi\circ\theta$ agrees with the regular
  representation $\L=(\pi_\AA\otimes\l)\circ \d\times (1\otimes M)$ on
  the generators $j_A(f)j_G(g)$ for $f\in \sa_c(G,\AA)$ and $g\in
  C_c(G)$.  Hence $\Pi\circ\theta = \L$ on all of $\cs(G,\AA)\times_\d
  G$ by linearity, continuity, and density.
\end{proof}

\section{Semidirect-product bundles}
\label{Semidirect-product bundles}

To
prove our main theorem in Section~\ref{canonical}, we are going to
need to build a Fell bundle over groupoid arising as a semidirect
product.  In this section, we give the construction of this
\emph{semidirect-product Fell bundle}.  We will investigate the
structure of the corresponding Fell bundle $C^{*}$-algebra in
Section~\ref{action crossed product}.

To begin, let $\GG$ be a locally compact Hausdorff groupoid with Haar
system $\set{\lambda^{u}}_{u\in \GG^{(0)}}$, and let $G$ be a 
second countable locally compact group.  
An \emph{action} of $G$ on $\GG$ is a homomorphism
$\b:G\to\aut\GG$ such that
$(x,t)\mapsto \b_t(x)$ is continuous from $\GG\times G$ to $\GG$.
(Note that
automorphisms of a groupoid do not necessarily fix the unit space
pointwise.) 
Given an action $\beta$ of $G$ on $\GG$,
the \emph{semidirect-product groupoid} $\GG\times_\b G$ comprises the
Cartesian product $\GG\times G$ with multiplication
\[
(x,t)(y,s)=(x\b_t(y),ts)
\]
whenever $s(x)=\beta_{t}(r(y))$
and inverse $(x,t)\inv=(\b_{t\inv}(x\inv),t\inv)$ 
(\cite[Definition~I.1.7]{ren:approach}). 
Note that we have
$(\GG\times_\b G)^0=\GG^0\times \{e\}$, with
\[
  r(x,t)=(r(x),e)
\quad\text{and}\quad
  s(x,t)=(\b_t\inv(s(x)),e).
\]
Also note that $C_c(\GG)\odot C_c(G)$ is inductive-limit dense in
$C_c(\GG\times_\b G)$.

Now suppose $p:\BB\to \GG$ is a separable
Fell bundle over $\GG$.  An \emph{action} of $G$ on $\BB$ is a homomorphism
$\a:G\to \aut\BB$ such that
$(b,t)\mapsto \a_t(b)$
is continuous from  $\BB\times G\to\BB$,
together with an \emph{associated action} $\b$ of $G$ on $\GG$
such that
$p\bigl(\alpha_{t}(b)\bigr)= \beta_{t}\bigl(p(b)\bigr)$
for all $t\in G$ and $b\in\BB$. 

\begin{remark}
  \label{rem-alpha-bar}
  The compatibility of $\alpha$ and $\beta$ allows us to write down,
  for each $t\in G$, an automorphism $\bar\alpha_t$ of
  $\sa_{c}(\GG,\BB)$ given by
  \begin{equation}
    \label{eq:13}
    \bar\alpha_t(f)(x)=\alpha_{t}\bigl(f(\beta_{t}^{-1}(x)\bigr).
  \end{equation}
  Since $\bar\alpha_{t}$ is clearly continuous from the inductive
  limit topology to the norm topology, it follows from
  Proposition~\ref{extends} that $\bar\alpha_{t}$ extends to an
  automorphism of $\cs(\GG,\BB)$.  Similarly, $t\mapsto
  \bar\alpha_{t}(f)$ is continuous from $G$ into $\cs(\GG,\BB)$, so
  we obtain an action $\bar\alpha$ of $G$ on $\cs(\GG,\BB)$.
\end{remark}

\begin{prop}
  \label{prop-fb-needed}
Let $\alpha$ be an action of $G$ on a Fell bundle $p: \AA\to \GG$,
with associated action $\beta$ of $G$ on $\GG$.  
Then the Banach bundle 
$q:\BB\times_\a G\to \GG\times_\b G$ with total space 
$\BB\times G$ and bundle projection  $q(b,t) = (p(b),t)$
becomes a Fell bundle when equipped 
 with the multiplication given by
  \[
  (b_x,t)(c_y,s)=(b_x\a_t(c_y),ts)
\quad\text{whenever $s(x)=r(\b_t(y))$}
  \]
  and the involution given by
  \[
  (b_x,t)^*=(\a_{t\inv}(b_x)^*,t\inv).
  \]
\end{prop}

We refer to a Fell bundle which arises from a group action as in Proposition~\ref{prop-fb-needed} as a
\emph{semidirect-product Fell bundle}.  

\begin{proof}[Sketch of Proof]
  For convenience, we'll write $C_{(x,t)}$ for the fibre of $\btaa$
  over $(x,t)$.  Verifying the axioms that $\btaa$ is a Fell bundle is
  routine with the possible exception of seeing that $C_{(x,t)}$ is a
  $C_{(r(x),e)}$\,--\,$C_{(\beta_{t}^{-1}(s(x)),e)}$-imprimitivity
  bimodule with respect to the operations inherited from $\btaa$.
  However, $C_{(x,t)}$ is naturally identified with $B_{x}$, and the
  latter is given to be a $B_{r(x)}$\,--\,$B_{s(x)}$-imprimitivity
  bimodule with respect to the operations inherited from $\BB$.
  Furthermore, $\alpha_{t}$ restricts to a $C^{*}$-algebra isomorphism
  of $B_{\beta_{t}^{-1}(s(x))}$ onto $B_{s(x)}$.  Therefore $B_{x}$ is
  naturally a $B_{r(x)}$\,--\,$B_{\beta_{t}^{-1}(s(x))}$-imprimitivity
  bimodule.  The right action is given by $x\cdot b=x\alpha_{t}(b)$
  and the right inner product is given by
  \[
  \langle x,y\rangle_{\ripsqueeze\ipscriptstyle
    B_{\beta_{t}^{-1}(s(x))}} =\alpha_{t}^{-1}\bigl(\langle x\ipcomma
  y \rangle_{\ripsqueeze\ipscriptstyle{B_{s(x)}}}\bigr).
  \]
  Now it is a simple matter to see that the given operations in $\btaa$
  induce the same structure on $C_{(x,t)}$ as does
the identification of $C_{(x,t)}$ with $B_x$.
\end{proof}

In order to have a Haar system on a semidirect-product groupoid
$\GG\times_\beta G$,
we will  need $\beta$ to be compatible with
the Haar system on $\GG$ in the following sense.

\begin{defn}\label{invt-def}
  An action $\b:G\to\aut\GG$ is \emph{invariant} if for all
  $u\in\go$, $f\in C_c(\GG)$, and $t\in G$ we have
  \[
  \int_\GG f(\b_t(y))\,d\l^u(y)=\int_\GG f(y)\,d\l^{\b_t(u)}(y),
  \]
  i.e., $\b_t$ transforms the measure on $r\inv(u)$ to the
  measure on $r\inv(\b_t(u))$.
If $\a:G\to\aut\BB$ is an action on a Fell
  bundle $\BB\to \GG$ with associated action $\b:G\to \aut\GG$, we say
  $\a$ is \emph{invariant} if $\b$ is.
\end{defn}

\begin{prop}\label{lem-haar-system-gltgg}
  Let $\b:G\to\aut\GG$ be an invariant action on a groupoid $\GG$ with
  Haar system $\{ \l^u \}_{u\in \GG^{(0)}}$. Then
  \[
  d\l^{(u,e)}(y,s)=d\l^u(y)\,ds
  \]
  is a Haar system on $\GG\times_\b G$.
\end{prop}

\begin{proof}
  The left-invariance property we need is that for $h\in
  C_c(\GG\times_\b G)$ and $(x,t)\in \GG\times_\b G$ we have
  \[
  \int_{\GG\times G}h\bigl((x,t)(y,s)\bigr)\,d\l^{s(x,t)}(y,s)
  =\int_{\GG\times G}h(y,s)\,d\l^{r(x,t)}(y,s),
  \]
  and it suffices to take $h=f\otimes g$, where $f\in G_c(\GG)$ and
  $g\in C_c(G)$.  Fix $x\in \GG$ with $s(x)=v$ and $r(x)=u$.  In the
  left-hand integral we must have
  \[
  (r(y),e)=r(y,s)=s(x,t)=(\b_t\inv(s(x)),e)=(\b_t\inv(v),e),
  \]
  and in the right-hand integral we must have
  \[
  (r(y),e)=(r(x),e)=(u,e).
  \]
  Since
  \[
  (x,t)(y,s)=(x\b_t(y),ts),
  \]
  we must show that
  \[
  \int_\GG\int_Gf(x\b_t(y))g(ts)\,ds\,d\l^{\b_t\inv(v)}(y)
  =\int_\GG\int_Gf(y)g(s)\,ds\,d\l^u(y).
  \]
  We have
  \begin{align*}
    \int_\GG\int_Gf(x\b_t(y))g(ts)\,ds\,d\l^{\b_t\inv(v)}(y)
    &=\int_\GG f(x\b_t(y))\int_G g(ts)\,ds\,d\l^{\b_t\inv(v)}(y)
    \\&=\int_\GG f(x\b_t(y))\,d\l^{\b_t\inv(v)}(y)\int_Gg(s)\,ds
  \end{align*}
  and similarly
  \begin{align*}
    \int_\GG\int_Gf(y)g(s)\,ds\,d\l^u(y) =\int_\GG
    f(y)\,d\l^u(y)\int_G g(s)\,ds,
  \end{align*}
  so it remains to verify
  \[
  \int_\GG f(x\b_t(y))\,d\l^{\b_t\inv(v)}(y)=\int_\GG f(y)\,d\l^u(y).
  \]
  But invariance of the action $\b$ gives
  \[
  \int_\GG f(x\b_t(y))\,d\l^{\b_t\inv(v)}(y)=\int_\GG f(xy)\,d\l^v(y),
  \]
  which equals $\int_\GG f(y)\,d\l^u(y)$ because $\l$ is a Haar
  system.
\end{proof}

For reference, we record the formula for convolution in
$C_c(\GG\times_\b G)$:
\begin{align*}
  &(h*k)(x,t) =\int_\GG\int_G h(y,s)k\bigl(\b_s\inv(y\inv x),s\inv
  t\bigr)\,ds\,d\l^{r(x)}(y).
\end{align*}
Thus in $\sa_c(\GG\times_\beta G,\BB\times_\a G)$ the convolution is
given by
\begin{align*}
  (h*k)(x,t) &=\int_\GG\int_G h(y,s)k\bigl(\b_s\inv(y\inv
  x),s\inv t\bigr)\,ds\,d\l^{r(x)}(y) \\
  &=\int_\GG\int_G (h_1(y,s),s)\bigl(k_1\bigl(\b_s\inv(y\inv x),s\inv
  t\bigr),s\inv
  t\bigr)\,ds\,d\l^{r(x)}(y)  \\
  &=\int_\GG\int_G \bigl( h_1(y,s)\a_s\bigl(k_1\bigl(\b_s\inv(y\inv
  x),s\inv t\bigr)\bigr) ,t \bigr) \,ds\,d\l^{r(x)}(y).
\end{align*}
As with product bundles (see Section~\ref{product bundles}), every
section $h\in\sa_c(\GG\times_\b G,\BB\times_\a G)$ is of the form
\[
h(x,t)= (h_1(x,t), t),
\]
where $h_1\in C_c(\GG\times_\b G,\BB)$ satisfies $h_1(x,t) \in B_x$.
So in particular
\begin{equation}\label{hkxt}
  (h*k)_1(x,t)
  =\int_\GG\int_G
  h_1(y,s)\a_s\bigl(k_1\bigl(\b_s\inv(y\inv x),s\inv t\bigr)\bigr)
  \,ds\,d\l^{r(x)}(y).
\end{equation}
The involution in $\sa_c(\GG\times_\beta G,\BB\times_\a G)$ is given
by
\begin{align*}
  h^*(x,t) &=h\bigl((x,t)\inv\bigr)^*
  =h\bigl(\b_t\inv(x\inv),t\inv\bigr)^*
  =\bigl(h_1\bigl(\b_t\inv(x\inv),t\inv\bigr),t\inv\bigr)^*
  \\&=\bigl(\a_t\bigl(h_1\bigl(\b_t\inv(x\inv),t\inv\bigr)^*\bigr),t\bigr),
\end{align*}
so in particular
\[
h^*_1(x,t) =\a_t\bigl(h_1\bigl(\b_t\inv(x\inv),t\inv\bigr)^*\bigr).
\]


\section{Action crossed product}
\label{action crossed product}

We now relate the $C^*$-algebra of a semidirect-product bundle to the
crossed product.

\begin{thm}
  \label{action crossed product isomorphism}
  Let $p:\BB\to\GG$ be a separable Fell bundle over a locally compact
  Hausdorff groupoid with Haar system $\{ \lambda^u \}_{u\in G^{(0)}}$, 
and let $\a:G\to\aut \BB$ be an action of a
  second countable locally compact group $G$ on $\BB$ with an invariant associated
  action $\b$ of $G$ on $\GG$. Let $q:\BB\times_\a G\to
  \GG\times_\b G$ denote the associated semidirect-product Fell bundle
  over the semidirect-product groupoid as defined in
  Section~\ref{Semidirect-product bundles}, and let $\bar\a:G\to\aut
  \cs(\GG,\BB)$ denote the concomitant action described in
  Remark~\ref{rem-alpha-bar}.  Then there is a unique isomorphism
  \[
  \s:\cs(\GG,\BB)\rtimes_{\bar\a} G\longrightarrow \cs(\GG\times_\b
  G,\BB\times_\a G)
  \]
  such that if $f\in \sa_c(\GG,\BB)$ and $g\in C_c(G)$ then
  $\s(i_\BB(f)i_G(g))$ is the continuous compactly supported section
  of $\BB\times_\a G$ given by
  \begin{equation}\label{s-eq1}
  \s\bigl(i_\BB(f)i_G(g)\bigr)(x,t)=\bigl(f(x)g(t)\D(t)^{\half},t\bigr).
  \end{equation}
\end{thm}

\begin{proof}
  Uniqueness is immediate from density.  For existence,
we will obtain $\s$ as the integrated form of a covariant
homomorphism $(\s_\BB,\s_G)$ of $(\cs(\GG,\BB),G,\bar\a)$ 
  into $M(\cs(\GG\times_\b G,\BB\times_\a G))$ such that
\begin{equation}\label{s-eq2}
  \s_\BB(f)\s_G(g) = f\boxtimes(\D^{\half}g) \in \sa_c(\GG\times_\beta
  G,\BB\times_\a G)
\end{equation}
for $f\in\sa_c(\GG,\BB)$ and $g\in C_c(G)$.  
It will follow that $\s$ maps 
$\cs(\GG,\BB)\rtimes_{\bar\a}G$ into $\cs(\GG\times_\b G,\BB\times_\a G)$,
satisfies~\eqref{s-eq1}, and is surjective because
the sections in~\eqref{s-eq2} have
inductive-limit-dense span in $\sa_c(\GG\times_\beta G,\BB\times_\a G)$.

  To define $\s_{\BB}$, we will appeal to
  Proposition~\ref{prop-Fell-multipliers}, viewing $\cs(\GG\times_\b
  G,\BB\times_{\a}G)$ as a right Hilbert module over itself, with
  dense subspace $\sa_c(\GG\times_\beta G,\BB\times_\a G)$.  For $f\in
  \sa_{c}(\GG,\BB)$ we define a linear operator $\s_\BB(f)$ on
  $\sa_c(\GG\times_\beta G,\BB\times_\a G)$ by
  \begin{align*}
    \bigl(\s_\BB(f)h\bigr)(y,t)&=\int_\GG \bigl(f(x)h_1(x\inv
    y,t),t\bigr)\,d\l^{r(y)}(x).
  \end{align*}
  Seeing that $\s_{\BB}: \sa_c(\GG,\BB)\to\lin(\sa_c(\GG\times_\beta
  G,\BB\times_\a G))$ is an algebra homomorphism is straightforward:
  for $f,g\in\sa_c(\GG,\BB)$ and $h\in \sa_c(\GG\times_\beta
  G,\BB\times_\a G)$ we have
  \begin{align*}
    \bigl(\s_\BB(f)\s_\BB(g)h\bigr)(y,t) &=\int_\GG\bigl(
    f(x)\bigl(\s_\BB(g)h\bigr)_1(x\inv y,t)
    ,t\bigr)\,d\l^{r(y)}(x) \notag \\
    &=\int_\GG\int_\GG\bigl( f(x)g(z)h_1(z\inv x\inv y,t) ,t\bigr)
    \,d\l^{s(x)}(z) \,d\l^{r(y)}(x), \intertext{which, after using
      Fubini and sending $z\mapsto x\inv z$, is} &=\int_\GG\int_\GG
    \bigl(f(x)g(x\inv z) h_1(z\inv y,t),t\bigr)
    \,d\l^{r(y)}(x) \,d\l^{r(y)}(z)\\
    &=\int_\GG\bigl(f*g(z)h_1(z\inv y,t),t\bigr)\,d\l^{r(y)}(z) \\
    &=\bigl(\s_\BB(f*g)h\bigr)(y,t).
  \end{align*}
  Thus, it remains to verify that $\s_{\BB}$ satisfies (i), (ii) and
  (iii) of Proposition~\ref{prop-Fell-multipliers}.  

To check (i), we
  compute as follows.  For $h,k\in \sa_c(\GG\times_\beta
  G,\BB\times_\a G)$, we have
  \begin{align*}
    \brip{\scriptstyle1}<\s&_\BB(f)h,k>(x,t)
    = \bigl( (\s_\BB(f)h)^**k)_1(x,t)\\
    &= \int_\GG\int_G (\s_\BB(f)h)^*_1(y,s) \a_s\bigl(
    k_1(\b_s\inv(y\inv x),s\inv t)\bigr)
    \,ds\,d\l^{r(x)}(y)\\
    &= \int_\GG\int_G \a_s\bigl(
    (\s_\BB(f)h)_1(\b_s\inv(y\inv),s\inv)\bigr)^* \a_s\bigl(
    k_1(\b_s\inv(y\inv x),s\inv t)\bigr)
    \,ds\,d\l^{r(x)}(y)\\
    &= \int_\GG\int_G\int_\GG \a_s\bigl(
    f(z)h_1(z\inv\b_s\inv(y\inv),s\inv)\bigr)^*
    \a_s\bigl( k_1(\b_s\inv(y\inv x),s\inv t)\bigr)\\
    &\hspace{3in} \,d\l^{r(\b_s\inv(y\inv))}(z)\, ds\,d\l^{r(x)}(y)\\
    &= \int_\GG\int_G\int_\GG \a_s\Bigl(
    h_1(z\inv\b_s\inv(y\inv),s\inv)^*f(z)^*
    k_1(\b_s\inv(y\inv x),s\inv t)\Bigr)\\
    &\hspace{3in}\,d\l^{r(\b_s\inv(y\inv))}(z)\, ds\,d\l^{r(x)}(y)\\
    \intertext{which, after $z\mapsto \b_s\inv(y\inv)z$ for fixed $y$,
      is} &= \int_\GG\int_G\int_\GG \a_s\Bigl(
    h_1(z\inv,s\inv)^*f(\b_s\inv(y\inv)z)^*
    k_1(\b_s\inv(y\inv x),s\inv t)\Bigr)\\
    &\hspace{3in}\,d\l^{s(\b_s\inv(y\inv))}(z)\, ds\,d\l^{r(x)}(y)\\
    \intertext{which, by invariance of the action $\b$ (in the variable $z$), is} &=
    \int_\GG\int_G\int_\GG \a_s\Bigl(
    h_1(\b_s\inv(z\inv),s\inv)^*f(\b_s\inv(y\inv z))^*
    k_1(\b_s\inv(y\inv x),s\inv t)\Bigr)\\
    &\hspace{3in}\,\l^{s(y\inv)}(z)\, ds\,d\l^{r(x)}(y)\\
    \intertext{which, by Fubini, is} &= \int_\GG\int_G\int_\GG
    \a_s\Bigl( h_1(\b_s\inv(z\inv),s\inv)^*f(\b_s\inv(y\inv z))^*
    k_1(\b_s\inv(y\inv x),s\inv t)\Bigr)\\
    &\hspace{3in}\,d\l^{r(x)}(y)\, ds\,\l^{r(x)}(z)\\
    \intertext{which, after $y\mapsto zy$ for fixed $z$, is} &=
    \int_\GG\int_G\int_\GG \a_s\Bigl(
    h_1(\b_s\inv(z\inv),s\inv)^*f(\b_s\inv(y\inv))^*
    k_1(\b_s\inv(y\inv z\inv x),s\inv t)\Bigr)\\
    &\hspace{3in}\,d\l^{s(z)}(y)\, ds\,\l^{r(x)}(z)\\
    \intertext{which, by invariance of $\b$ (in $y$), is} &=
    \int_\GG\int_G\int_\GG \a_s\Bigl(
    h_1(\b_s\inv(z\inv),s\inv)^*f(y\inv)^*
    k_1(y\inv\b_s\inv(z\inv x),s\inv t)\Bigr)\\
    &\hspace{3in}\,d\l^{s(\b_s\inv(z))}(y)\, ds\,\l^{r(x)}(z)\\
    &= \int_\GG\int_G\int_\GG \a_s\bigl(
    h_1(\b_s\inv(z\inv),s\inv)^*\bigr)
    \a_s\bigl( f^*(y)k_1(y\inv \b_s\inv(z\inv x),s\inv t)\bigr)\\
    &\hspace{3in}\,d\l^{r(\b_s\inv(z\inv x))}(y)\,ds\,d\l^{r(x)}(z)\\
    &= \int_\GG\int_G h_1^*(z,s)\a_s\bigl(
    (\s_\BB(f^*)k)_1(\b_s\inv(z\inv x),s\inv t)\bigr)
    \,ds\,d\l^{r(x)}(z)\\
    &= \bigl( h^**(\s_\BB(f^*)k)\bigr)_1(x,t)\\
    &= \brip{\scriptstyle1}< h,\s_\BB(f^*)k>(x,t).
  \end{align*}

  To check the continuity condition~(ii) of
  Proposition~\ref{prop-Fell-multipliers}, it suffices to show that if
  $L\subseteq \GG$ is compact and $f_i\to 0$ uniformly in
  $\G_L(\GG,\BB)$, then for each $h,k\in \sa_c(\GG\times_\b
  G,\BB\times_\a G)$ there exists a compact set $K\subseteq
  \GG\times_\b G$ such that $\brip< \s_\BB(f_i)h,k>\to 0$ uniformly in
  $\G_K(\GG\times_\b G,\BB\times_\a G)$.  Using continuity of the
  action of $G$ on $\GG$, it is routine to verify that for any such
  $h$ and $k$ there exists a compact set $K$ such that $\supp \brip<
  \s_\BB(f_i)h,k> \subseteq K$ for every $i$.  Then, to verify uniform
  convergence, we notice that for each $i$,
  \begin{equation*}
    \|\brip< \s_{\BB}(f_i)h,k>\|_{\infty}\le M
    \|f_i\|_{\infty}\|h\|_{\infty} \|k\|_{\infty},
  \end{equation*}
  where $M=\sup_{u\in \GG^{(0)}}\lambda^{(e,u)}(K)$.

  For the nondegeneracy condition~(iii) of
  Proposition~\ref{prop-Fell-multipliers}, note that if
  $f,g\in\sa_c(\GG,\BB)$ and $h\in C_c(G)$, then
  \[
  \s_\BB(f)(g\boxtimes h)=(f*g)\boxtimes h,
  \] where $g\boxtimes h\in\sa_c(\GG\times_\beta G,\BB\times_\a G)$ is
  defined by $(g\boxtimes h)(x,t) = (g(x)h(t),t)$.
  Letting $f$ run through an approximate identity $\{f_i\}$ for
  $\sa_c(\GG,\BB)$ in the inductive limit topology (see
  \cite[Proposition~6.10]{mw:fell}), we have $f_i*g\to g$, hence
  $(f_i*g)\boxtimes h\to g\boxtimes h$, both nets converging in the
  inductive limit topology. Since such sections $g\boxtimes h$ have
  dense span in $\sa_c(\GG\times_\beta G,\BB\times_\a G)$, hence in
  $\cs(\GG\times_\b G,\BB\times_\a G)$, nondegeneracy follows.

  Now we conclude from Proposition~\ref{prop-Fell-multipliers} that
  $\s_{\BB}$ extends to a nondegenerate $*$-ho\-mo\-mor\-phism of
  $\cs(\GG,\BB)$ into $M(\cs(\GG\times_\b G,\BB\times_\a G))$, as
  required.


  We now turn to $\s_G$.  Fix $s\in G$, and for each $h\in
  \sa_c(\GG\times_\beta G,\BB\times_\a G)$, define $\s_G(s)h \in
  \sa_c(\GG\times_\beta G,\BB\times_\a G)$ by
  \begin{align*}
    \bigl(\s_G(s)h\bigr)(x,t)
    &=\bigl(\a_s\bigl(h_1\bigl(\b_s\inv(x),s\inv
    t\bigr)\bigr)\D(s)^{\half},t\bigr).
  \end{align*}
  Then for $h,k\in \sa_c(\GG\times_\beta G,\BB\times_\a G)$ we have
  \begin{align*}
    \brip{\scriptstyle 1}< \s_G&(s)h, \s_G(s)k>(x,t)
    = \bigl( (\s_G(s)h)^* * (\s_G(s)k) \bigr)_1(x,t)\\
    &= \int_\GG\int_G (\s_G(s)h)^*_1(y,r) \a_r\bigl(
    (\s_G(s)k)_1(\b_r\inv(y\inv x),
    r\inv t)\bigr) \, dr\, d\l^{r(x)}(y)\\
    &= \int_\GG\int_G
    \a_r\bigl( (\s_G(s)h)_1(\b_r\inv(y\inv),r\inv)^*\bigr) \\
    &\hspace{1.5in} \a_r\bigl( (\s_G(s)k)_1(\b_r\inv(y\inv x), r\inv
    t)\bigr) \, dr\,
    d\l^{r(x)}(y)\\
    &= \int_\GG\int_G \a_r\bigl(
    \a_s(h_1(\b_s\inv(\b_r\inv(y\inv)),s\inv r\inv)^*\D(s)^{\half}
    \bigr)\\
    &\hspace{1.5in} \a_r\bigl( \a_s(k_1(\b_s\inv(\b_r\inv(y\inv x)),
    s\inv r\inv t)\D(s)^{\half}
    \bigr) \, dr\, d\l^{r(x)}(y)\\
    &= \int_\GG\int_G
    \a_{rs}\bigl( h_1(\b_{rs}\inv(y\inv),(rs)\inv)^*\bigr) \\
    &\hspace{1.5in} \a_{rs}\bigl( k_1(\b_{rs}\inv(y\inv x), (rs)\inv
    t)\bigr) \,
    \D(s)dr\, d\l^{r(x)}(y)\\
    \intertext{which, after $r\mapsto rs^{-1}$, is} &= \int_\GG\int_G
    \a_{r}\bigl( h_1(\b_{r}\inv(y\inv),r\inv)^*\bigr)
    \a_{r}\bigl( k_1(\b_{r}\inv(y\inv x), r\inv t)\bigr) \, dr\, d\l^{r(x)}(y)\\
    &=\int_\GG\int_G h_1^*(y,r) \a_r\bigl( k_1(\b_r\inv(y\inv x),r\inv
    t)\bigr)\, dr\,
    d\l^{r(x)}(y)\\
    &= ( h^**k )_1(x,t) = \rip{\scriptstyle1}< h,k >(x,t).
  \end{align*}
  Since we clearly have $\s_{G}(s)\s_{G}(t)=\s_{G}(st)$
    and $\s_{G}(e)$ is the identity, it follows that $\s_{G}(s)$
    defines a unitary in $M(\cs(\GG\times_\b G,\BB\times_\a G))$.

  To see that
the resulting homomorphism
$\s_G:G\to M(\cs(\GG\times_\b G,\BB\times_\a G))$
is strictly continuous, it suffices (by
  \cite[Corollary~C.8]{tfb}) to show that if $s_i\to e$ in $G$ and
  $h\in\sa_c(\GG\times_\beta G,\BB\times_\a G)$ then $\s_G(s_i)h\to h$
  in the inductive limit topology. Without loss of generality all the
  $s_i$'s are contained in some compact neighborhood $V$ of
  $e$. Choose compact sets $K\subseteq \GG$ and $L\subseteq G$ such that
  $\supp h\subseteq K\times L$.  Then for each $i$ we have
  \[
  \supp \s_G(s_i)h\subseteq \b_{V\inv}(K)\times V\inv L,
  \]
  which is compact by continuity of the action $\b$. The uniform
  continuity of $h$ and continuity of the actions $\a$ and $\b$
  guarantee that
  \[
  \lim_i\a_{s_i}\bigl(h_1(\b_{s_i\inv}(x),s_i\inv t\bigr)\bigr) =
  h_1(x,t)
  \]
  uniformly in $(x,t)$, so $\s_G(s_i)h\to h$ uniformly. Thus
  $\s_G(s_i)h\to h$ in the inductive limit topology.

  Now we verify that the pair $(\s_\BB,\s_G)$ is covariant
for $(\cs(\GG,\BB),G,\bar\a)$.  If $f\in
  \sa_c(\GG,\BB)$ and $s\in G$, then for each $h\in
  \sa_c(\GG\times_\beta G,\BB\times_\alpha G)$ and $(y,t)\in
  \GG\times_\beta G$, we have
  \begin{align*}
    \bigl(\s_G(s)&\s_\BB(f)h\bigr)_1(y,t)
    = \a_s\bigl((\s_\BB(f)h)_1(\b_s\inv(y),s\inv t)\bigr)\D(s)^{\half}\\
    &= \int_\GG \a_s\bigl( f(x)h_1(x\inv\b_s\inv(y),s\inv
    t)\bigr)\D(s)^{\half}
    \,d\l^{r(\b_s\inv(y))}(x)\\
    \intertext{which, by invariance of $\b$, is} &= \int_\GG
    \a_s\bigl( f(\b_s\inv(x))h_1(\b_s\inv(x\inv y),s\inv
    t)\bigr)\D(s)^{\half}
    \,d\l^{r(y)}(x)\\
    &= \int_\GG \bar\a_s(f)(x)\bigl(\s_G(s)h\bigr)_1(x\inv y,t)
    \,d\l^{r(y)}(x)\\
    &= \bigl(\s_\BB(\bar\a_s(f))\s_G(s)h\bigr)_1(y,t).
  \end{align*}

Next we verify~\eqref{s-eq2}:
for
  $h\in\sa_c(\GG\times_\beta G,\BB\times_\a G)$ and $(y,s)\in
  \GG\times_\beta G$, we have
  \begin{align*}
    \bigl(\s_\BB(f)\s_G(g)&h\bigr)(y,s) =\Bigl( \int_\GG f(x)
    \bigl(\s_G(g)h\bigr)_1(x\inv y,s) \,d\l^{r(y)}(x) ,s\Bigr)
    \\
    &=\Bigl( \int_\GG \int_G f(x) g(t) \a_t\Bigl(
    h_1\bigl(\b_t\inv(x\inv y),t\inv s) \Bigr)
    \D(t)^{\half}\,dt \,d\l^{r(y)}(x) ,s\Bigr) \\
    &= \int_\GG\int_G \bigl( f(x)g(t)\a_t\bigl(h_1(\b_t\inv(x\inv
    y),t\inv
    s)\bigr),s\bigr)\D(t)^{\half} \,dt\,d\l^{r(y)}(x) \\
    & = \int_\GG\int_G
    \bigl(f(x)g(t)\D(t)^{\half},t\bigr)\bigl(h_1(\b_t\inv(x\inv
    y),t\inv
    s),t\inv s\bigr) \,dt\,d\l^{r(y)}(x) \\
    & = \int_\GG\int_G \bigl(f\boxtimes(\D^{\half}g)\bigr)(x,t)
    h(\b_t\inv(x\inv y),t\inv
    s) \,dt\,d\l^{r(y)}(x) \\
    & = \bigl(\bigl(f\boxtimes(\D^{\half}g)\bigr)*h\bigr)(y,s).
  \end{align*}

As outlined at the start of the proof, it follows from the above
that the integrated form $\s=\s_\BB\rtimes \s_G$  maps
  $\cs(\GG,\BB)\rtimes_{\bar\a} G$ (into and) onto $\cs(\GG\times_\b
  G,\BB\times_\a G)$. 
  To show injectivity of $\sigma$, it suffices to find a left inverse.
  We will begin by constructing a $*$-homomorphism
  $\tau:\sa_c(\GG\times_\beta G,\BB\times_\a G)\to
  C_c(G,\sa_c(\GG,\BB))$ which is continuous for the inductive limit
  topologies on each algebra, where (of course) $\sa_c(\GG,\BB)$ is
  also given the inductive limit topology.  Then, the composition
  \[
  \sa_c(\GG\times_\b G,\BB\times_\a G) \xrightarrow{\tau}
  C_c(G,\sa_c(\GG,\BB)) \xrightarrow{} C_c(G,C^*(\GG,\BB))
  \xrightarrow{} C^*(\GG,\BB)\times_{\bar\a}G
  \]
  will be continuous from the inductive limit topology to the
  $C^*$-norm topology, and hence, by Proposition~\ref{extends}, will
  extend to a homomorphism, which we will also denote by $\tau$, of
  $C^*(\GG\times_\b G,\BB\times_\a G)$ into
  $C^*(\GG,\BB)\rtimes_{\bar\a}G$.  Finally, we will check
  that $\tau\circ\sigma=\id$ on generators, and this will suffice.

  For $h\in \sa_c(\GG\times_\beta G,\BB\times_\a G)$ and $t\in G$, it
  is clear that the rule
  \[
  x\mapsto h_1(x,t)\D(t)^{-\half}
  \]
  defines an element $\tau(h)(t)$ of $\sa_c(\GG,\BB)$.  The discussion
  in \cite[II.15.19]{fd1} shows that the map $t\mapsto h_1(\cdot,t)$
  from $G$ into $\sa_c(\GG,\BB)$ is inductive-limit continuous, and it
  follows that $t\mapsto \tau(h)(t)$ defines an inductive-limit
  continuous map $\tau(h)$ from $G$ to $\sa_c(\GG,\BB)$.  Since
  $\tau(h)$ obviously has compact support, we therefore have
  $\tau(h)\in C_c(G,\sa_c(\GG,\BB))$, with
  \[
  \tau(h)(t)(x) = h_1(x,t)\D(t)^{-\half} \quad\text{for $t\in G$ and $
    x\in \GG$.}
  \]

  Now the rule $h\mapsto \tau(h)$ gives a map $\tau$ with domain
  $\sa_c(\GG\times_\b G,\BB\times_\a G)$ which is clearly linear.  To
  show that $\tau$ is continuous for the inductive limit topologies,
  it suffices to show that if $K\subseteq \GG$ and $L\subseteq G$ are
  compact and $\{ h_i\}$ is a net in $\G_{K\times L}(\BB\times_\a G)$
  converging uniformly to $0$, then $\tau(h_i)\to 0$ in the inductive
  limit topology of $C_c(G,\sa_c(\GG,\BB))$. Since
  $\supp\tau(h_i)\subseteq L$ for all $i$, it suffices to show that
  $\tau(h_i)\to 0$ uniformly. But this is obvious, since $h_i\to 0$
  uniformly.

  Next we show that $\tau$ is a $*$-homomorphism. For $h,k\in
  \sa_{c}(\GG\times_{\b}G,\BB\times_{\a}G)$ we can use the argument%
  \footnote{Lemma~1.108 of \cite{danacrossed} as stated does not apply
    to a section algebra $\sa_c(\GG,\BB)$ sitting inside a bundle
    $C^*$-algebra $\cs(\GG,\BB)$, but it is easy to see that the
    argument gives the conclusion we need here.}
of \cite[Lemma~1.108]{danacrossed} to conclude that $\tau(h)*\tau(k)$,
  which is \emph{a priori} an element of $C_{c}(G,\cs(\GG,\BB))$, lies
  in $C_{c}(G,\sa_{c}(\GG,\BB))$ and that we can pass ``evaluation at
  $x$'' through the integral in the second line of the following
  computation:
  \begin{align*}
    \bigl(\tau(h)&{}*\tau(k)\bigr)(t)(x)
    =\Bigl(\int_G\tau(h)(s)*\bar\a_s\bigl(\tau(k)(s\inv
    t)\bigr)\,ds\Bigr)(x)
    \\
    &=\int_G\bigl(\tau(h)(s)*\bar\a_s\bigl(\tau(k)(s\inv
    t)\bigr)\bigr)(x)\,ds \\
    &=\int_G\int_\GG\tau(h)(s)(y) \bar\a_s\bigl(\tau(k)(s\inv
    t)\bigr)(y\inv x)\,d\l^{r(x)}(y)\,ds
    \\
    &=\int_G\int_\GG h_1(y,s)\D(s)^{-\half} \a_s\bigl(\tau(k)(s\inv
    t)\bigl(\b_s\inv(y\inv
    x)\bigr)\bigr)\,d\l^{r(x)}(y)\,ds \\
    &=\int_G\int_\GG h_1(y,s)\D(s)^{-\half}
    \a_s\bigl(k_1\bigl(\b_s\inv(y\inv x),s\inv t\bigr)\bigr) \D(s\inv
    t)^{-\half}\,d\l^{r(x)}(y)\,ds
    \\
    &=(h*k)_1(x,t)\D(t)^{-\half} \\
    &=\tau(h*k)(t)(x),
  \end{align*}
  so $\tau$ is multiplicative.  For the involution, we have
  \begin{align*}
    \tau(h)^*(t)(x) &=\bar\a_t\bigl(\tau(h)(t\inv)^*\D(t\inv)\bigr)(x)
    \\&=\bar\a_t\bigl(\tau(h)(t\inv)^*\bigr)(x)\D(t\inv)
    \\&=\a_t\bigl(\tau(h)(t\inv)^*\bigl(\b_{t\inv}(x)\bigr)\bigr)\D(t\inv)
    \\&=\a_t\bigl(\tau(h)(t\inv)\bigl(\b_{t\inv}(x\inv)\bigr)^*\bigr)\D(t\inv)
    \\&=\a_t\bigl(h_1\bigl(\b_{t\inv}
    (x\inv),t\inv\bigr)\D(t)^{\half}\bigr)^*\D(t\inv)
    \\&=\a_t\bigl(h_1\bigl(\b_{t\inv}(x\inv),t\inv\bigr)\bigr)^*\D(t)^{-\half}
    \\&=(h^*)_1(x,t)\D(t)^{-\half} \\&=\tau(h^*)(t)(x).
  \end{align*}
  Finally, we check $\tau\circ\s=\id$ on generators of the form
  $i_\BB(f)i_G(g)$ for $f\in\sa_c(\GG,\BB)$ and $g\in C_c(G)$:
  \begin{align*}
    \tau\circ\s\bigl(i_\BB(f)i_G(g)\bigr)(t)(x) &=\tau\bigl(f\boxtimes
    (\D^{\half}g)\bigr)(t)(x) \\
    &=\bigl((f\boxtimes (\D^{\half}g)\bigr)_1(x,t)\D(t)^{-\half}
    \\
    &=f(x)g(t)
    \\
    &=\bigl(i_B(f)i_G(g)\bigr)(t)(x).\qed
  \end{align*}
  \renewcommand\qed{\relax}
\end{proof}

\section{The canonical surjection is injective}
\label{canonical}

The object of this section is to prove our main result:

\begin{thm}
  \label{Phi-iso}
  Let $\AA$ be a separable Fell bundle over a group $G$, and let
  $\delta$ be the associated coaction of $G$ on $C^*(G,\AA)$ as in
  Proposition~\ref{exists}.  Then the canonical surjection
  \[
  \Phi: \cs(G,\AA)\rtimes_\d G\rtimes_{\hat\d} G \to \cs(G,\AA)\otimes
  \KK(L^2(G))
  \]
  is an isomorphism; hence $\delta$ is maximal.
\end{thm}

To do this, we will factor $\Phi$ into three isomorphisms, each
involving the $C^{*}$-algebra of a Fell bundle over a groupoid.  These
isomorphisms will be presented in
Propositions~\ref{Theta}--\ref{Upsilon}.  We will use the following
notation for canonical maps related to the double-crossed product
$\cs(G,\AA)\rtimes_\d G\rtimes_{\hat\d} G$:
\begin{align*}
  k_\AA=i_{\cs(G,\AA)\rtimes_\d G}\circ
  j_{\cs(G,\AA)}:{}&\cs(G,\AA)\to M(\cs(G,\AA) \rtimes_\d
  G\rtimes_{\hat\d} G)\\
  k_{C(G)}=i_{\cs(G,\AA)\rtimes_\d G}\circ j_G:{}&C_0(G)\to
  M(\cs(G,\AA)
  \rtimes_\d G\rtimes_{\hat\d} G)\\
  k_G=i_G:{}&G\to M(\cs(A,\GG)\rtimes_\d G\rtimes_{\hat\d} G).
\end{align*}
Note that the double-crossed product is densely spanned by products of
the form
\[
k_\AA(f) k_{C(G)}(g) k_G(h) \quad\text{for $f\in \sa_c(G,\AA)$ and $
  g,h\in C_c(G)$.}
\]

Our first isomorphism involves an iterated product Fell bundle.  Let
$\AA\times_{\lt}G$ be the transformation Fell bundle over the
transformation groupoid $G\times_\lt G$, as in
Section~\ref{transformation bundles}.  The group $G$ acts on both
$G\times_\lt G$ and $\AA\times_\lt G$ by right translation in the
second co\"ordinate:
\[
(\id_G\times\rt)_r(s,t) = (s,tr\inv) \quad\text{and}\quad
(\id_\AA\times\rt)_r(a_s,t) = (a_s,tr\inv).
\]
Thus we get a semidirect-product Fell bundle $\altgg$;
for simplicity,
we will denote the corresponding semidirect-product groupoid
$(\gltg)\times_{\betart}G$ by~$\gltgg$.  

The action of $G$ on $\gltg$
is invariant in the sense of Definition~\ref{invt-def}, since for each
$(e,u)\in (\gltg)^0 = \{e\}\times G$, $f\in C_c(\gltg)$, and $r\in G$
we have
\begin{align*}
  \int_{\gltg} f\bigl((\betart)_r(s,t)\bigr)\,d\l^{(e,u)}(s,t)
  &= \int_G f((\betart)_r(s,s\inv u)\,ds\\
  &= \int_G f(s,s\inv u r\inv)\,ds\\
  &= \int_{\gltg} f(s,t)\,d\l^{(e,ur\inv)}(s,t)\\
  &= \int_{\gltg} f(s,t)\,d\l^{(\betart)_r(e,u)}.
\end{align*} 
Therefore
  Proposition~\ref{lem-haar-system-gltgg} gives a Haar system on $\gltgg$,
so we can form the Fell-bundle $C^{*}$-algebra $\cs(\gltgg,\altgg)$.

\begin{prop}
  \label{Theta}
  There is an isomorphism
  \[
  \Theta:\cs(G,\AA)\rtimes_\d G\rtimes_{\hat\d} G \to
  \cs(\gltgg,\altgg)
  \]
  such that, for $f\in\sa_c(G,\AA)$ and $g,h\in C_c(G)$, the image
  $\Theta(k_\AA(f)k_{C(G)}(g)k_G(h))$ is in $\sa_{c}(\gltgg,\altgg)$,
  with
  \begin{equation}\label{generators}
    \Theta\bigl(k_\AA(f)k_{C(G)}(g)k_G(h)\bigr)(r,s,t)
    =\bigl(
    f(r)g(s)h(t)\D(rt)^{\half}
    ,s,t\bigr).
  \end{equation}
\end{prop}

\begin{proof}
  Theorem~\ref{coaction crossed product isomorphism} gives an
  isomorphism
  \[
  \t:\cs(G,\AA)\rtimes_\d G\to C^*(\gltg,\altg)
  \]
  such that
  \[
  \t(j_{\csga}(f)j_G(g))=(\D^{\half}f)\boxtimes g
  \]
  for $f\in \sa_c(G,\AA)$ and $g\in C_c(G)$.  We want to parlay this
  into our isomorphism $\Theta$.  First, we verify that $\t$ is
  equivariant for the dual action of $G$ on $\csga\rtimes_\d G$ and
  the action $\alphabarrt$
coming from
 the action 
  of $G$ on $\AA\times_{\lt} G$.  Note that for $h\in
  \sa_{c}(\gltg,\altg)$,
  \begin{equation*}
    \alphabarrt_{s}(h)(t,r)=\bigl(h_{1}(t,rs),r\bigr).
  \end{equation*}
  Thus for $f\in\sa_c(G,\AA)$, $g\in C_c(G)$ and $s\in G$ we have
  \begin{align*}
    \alphabarrt_s\circ\t\bigl(j_{\csga}(f)j_G(g)\bigr)(t,r)
    &=\alphabarrt_s\bigl((\D^{\half}f)\boxtimes g\bigr)(t,r)
    \\&=\bigl(\D(t)^{\half}f(t)g(rs),t\bigr)
    \\&=\bigl(\D(t)^{\half}f(t)\rt_s(g)(r),t\bigr)
    \\&=\bigl((\D^{\half}f)\boxtimes \rt_s(g)\bigr)(t,r),
    \intertext{so that}
    \alphabarrt_s\circ\t\bigl(j_{\csga}(f)j_G(g)\bigr)
    &=\t\bigl(j_{\csga}(f)j_G(\rt_s(g))\bigr)
    \\&=\t\Bigl(\hat\d_s\bigl(j_{\csga}(f)j_G(g)\bigr)\Bigr)
    \\&=\t\circ\hat\d_s\bigl(j_{\csga}(f)j_G(g)\bigr).
  \end{align*}
  Therefore we have an isomorphism
  \[
  \t\rtimes G:\cs(G,\AA)\rtimes_\d G\rtimes_{\hat\d} G \to
  \cs(\gltg,\altg)\times_{\alphabarrt} G.
  \]
  Now, Theorem~\ref{action crossed product isomorphism} gives an
  isomorphism
  \[
  \s:\cs(\gltg,\altg)\rtimes_{\alphabarrt} G \to \cs(\gltgg,\altgg)
  \]
  taking a generator $i_{\cs(\gltg,\altg)}(k)i_G(h)$ for $k\in
  \sa_c(\gltg,\altg)$ and $h\in C_c(G)$ to the section of $\altgg$
  given by
  \[
  \s\bigl(i_{C^*(G\times_\lt G,\AA\times_{\lt}
    G)}(k)\id_G(h)\bigr)(r,s,t)
  =\bigl(k(r,s)h(t)\D(t)^{\half},s,t\bigr).
  \]

  We now define $\Theta$ to be $\sigma\circ(\theta\rtimes G)$, and
  it only remains to verify \eqref{generators}. We have
  \begin{align*}
    \Theta\bigl(k_\AA(f)k_{C(G)}(g)&k_G(h)\bigr) =\s\circ (\t\times
    G)\bigl( i_{\cs(A,\GG)\times_\d
      G}\bigl(j_{\cs(G,\AA)}(f)j_G(g)\bigr)i_G(h) \bigr)
    \\
    &=\s\bigl( (\t\times G)\bigl( i_{\cs(A,\GG)\times_\d
      G}\bigl(j_{\cs(G,\AA)}(f)j_G(g)\bigr)i_G(h) \bigr) \bigr) \\
    &=\s\bigl( i_{C^*(G\times_\lt G,\AA\times_{\lt} G)}\circ \t\bigl(
    j_{\cs(G,\AA)}(f)j_G(g)
    \bigr)i_G(h) \bigr) \\
    &=\s\bigl( i_{C^*(G\times_\lt G,\AA\times_{\lt}
      G)}\bigl((\D^{\half}f)\boxtimes g\bigr)i_G(h) \bigr),
  \end{align*}
  so
  \begin{align*}
    \Theta\bigl(k_\AA(f)k_{C(G)}(g)k_G(h)\bigr)(r,s,t) &=\bigl(
    \bigl((\D^{\half}f)\boxtimes g\bigr)(r,s)h(t)\D(t)^{\half}
    ,t\bigr) \\&=\bigl( f(r)g(s)h(t)\D(rt)^{\half} ,s,t\bigr).
    \qed
  \end{align*}
  \renewcommand\qed{\relax}
\end{proof}

For our second isomorphism, we let $\EE$ denote the equivalence
relation groupoid $G\times G$ on the set $G$, and we endow $\EE$ with
the Haar system $\lambda^{(s,s)} = \delta_s\times \lambda$, where
$\delta_s$ is the point mass at $s$, and $\lambda$ is Haar measure on
$G$.  We then form the Cartesian product Fell bundle $\AA\times\EE$
over the Cartesian product groupoid $G\times\EE$, in analogy with the
group case in Section~\ref{product bundles}.

\begin{prop}
  \label{Psi}
  There is an isomorphism
  \begin{equation*}
    \Psi:\cs(\gltgg,\altgg)\to \cs(G\times\EE,\AA\times\EE)
  \end{equation*}
  such that, for $f\in\sa_{c}(\gltgg,\altgg)$, the image $\Psi(f)$ is
  in $\sa_c(G\times \EE, \AA\times\EE)$, with
  \begin{equation}\label{Psi-eq}
    \Psi(f)(r,s,t) =  (f_{1}(r,r^{-1}s,s^{-1}rt),s,t).
  \end{equation}
\end{prop}

\begin{proof}
  First notice that the groupoids
  $\gltgg=(G\times_{\lt}G)\times_{\alphart}G$ and $G\times \EE$ are
  isomorphic via the homeomorphism $\psi:\gltgg\to G\times\EE$ given
  by $\psi(r,s,t)=(r,rs,st)$.  Furthermore, the homeomorphism $\Psi_0:
  \altgg\to \AA\times\EE$ given by $\Psi_{0}(a_{r},s,t)=(a_{r},rs,st)$
  is a bundle map which covers $\psi$ and is an isometric isomorphism
  on each fibre.
  Routine computations show that $\Psi_0$ also preserves the
  multiplication and involution.
  Hence we can define a
  $*$-isomorphism $\Psi:\sa_{c}(\gltgg,\altgg)\to
  \sa_{c}(G\times\EE,\AA\times\EE)$ by
  \[
  \Psi(f)(r,s,t) = \Psi_0(f(\psi\inv (r,s,t)) =
  (f_{1}(r,r^{-1}s,s^{-1}rt),s,t).
  \]
  Because $\Psi_0$ is a homeomorphism, $\Psi$ is homeomorphic for the
  inductive limit topologies; therefore $\Psi$ extends to an
  isomorphism of the bundle $C^*$-algebras which
  satisfies~\eqref{Psi-eq}.
\end{proof}

\begin{prop}
  \label{Upsilon}
  There is an isomorphism
  \begin{equation*}
    \Upsilon:\cs(G\times\EE,\AA\times\EE)\to \cs(G,\AA)\tensor
    \KK\bigl(L^{2}(G)\bigr) 
  \end{equation*}
  such that, for every faithful nondegenerate representation
  $\pi:\cs(G,\AA)\to B(\HH)$, $f\in \sa_{c}(G\times\EE,\AA\times\EE)$,
  and $\xi\in C_c(G,\HH)$, we have
  \begin{equation}
    \label{eq:9}
    \bigl((\pi\tensor\id)\circ \Upsilon(f)\xi\bigr)(s) =\int_{G}
    \int_{G} \pi_{0}\bigl(f_{1}(r,s,t)\bigr)\xi(t)\D(r)^{-\half} \,dr
    \,dt,
  \end{equation}
  where $\pi_{0} = \pi\circ\iota$
as in Lemma~\ref{bundle map}.
\end{prop}

The proposition depends on the following lemma, which may be of
general interest.  As above, $\AA\times\GG$ denotes the Cartesian
product bundle over the Cartesian product groupoid $G\times\GG$.

\begin{lem}\label{omega}
  Let $\GG$ be a second countable locally compact groupoid such that
  $C^*(\GG)$ is nuclear.  There exists an an isomorphism
  $\o:\cs(G\times\GG,\AA\times \GG)\to \cs(G,\AA)\otimes C^*(\GG)$
  such that
  \[
  \o(g\boxtimes h)=(\D^{-\half}g)\otimes h \quad\text{for
    $g\in\sa_c(G,\AA)$ and $h\in C_c(\GG)$,}
  \]
  where $g\boxtimes h\in \sa_c(G\times\GG,\AA\times \GG)$ is defined
  by $(g\boxtimes h)(s,x) = (g(s)h(x),x)$.
\end{lem}

\begin{proof}
  For $a_{t}\in A_{t}$, define a linear operator $\pan(a_t)$ on
  $\sa_{c}(G\times\GG,\AA\times\GG)$ by\footnote{Although this
    construction almost exactly parallels that in Lemma~\ref{key}, we
    need to insert a modular function here (compare with \eqref{iota})
    because, just as with (\ref{eq:5}) in the proof of
    Theorem~\ref{coaction crossed product isomorphism}, there is no
    modular function in the definition of the involution in the
    $*$-algebra associated to a Fell bundle over a groupoid.}
  \begin{equation*}
    \bigl(\pan(a_{t})h\bigr)_1(s,x)=a_{t}h_1(t^{-1}s,x)\D(t)^{\half}.
  \end{equation*}
  Then a computation shows that
  \begin{equation*}
    \brip< \pan(a_{t})h,k> = \brip< h,\pan(a_{t}^{*})k>
    \quad\text{ for $h,k\in \sa_c(G\times\GG, \AA\times\GG)$,}
  \end{equation*}
  where we are viewing $\cs(G\times\GG,\AA\times\GG)$ as a right
  Hilbert module over itself with dense subspace
  $\sa_c(G\times\GG,\AA\times\GG)$.  Just as in the proof of
  Lemma~\ref{key}, it follows that $\pan(a_{t})$ is bounded as an
  operator on $\cs(G\times\GG,\AA\times\GG)$ with adjoint
  $\pan(a_{t}^{*})$, and that the rule $a_t \mapsto \pan(a_t)$
  therefore extends to a $*$-homomorphism $\pan$ of $\AA$ into
  $M(\cs(G\times\GG,\AA\times\GG))$.  The proof that $\pan$ is
  nondegenerate and strictly continuous also closely parallels the
  proof in Lemma~\ref{key} and will be omitted.  Using
  Lemma~\ref{bundle map}, we get a nondegenerate homomorphism $\pa:
  \cs(G,\AA)\to M(\cs(G\times\GG,\AA\times\GG))$.

  Similarly, for $g\in C_{c}(\GG)$ we define an operator $\pg(g)$ on
  $\sa_{c}(G\times\GG,\AA\times\GG)$ by
  \begin{equation*}
    (\pg(g)h)_1(s,x)=\int_{\GG} g(y)h_1(s,y^{-1}x)\,d\lambda^{r(x)}(y).
  \end{equation*}
  Another computation shows that
  \begin{equation*}
    \rip< \pg(g)h,k> =\brip <h, \pg(g^{*})k>
    \quad\text{ for $h,k\in \sa_c(G\times\GG, \AA\times\GG)$.}
  \end{equation*}
  Thus condition~(i) of Proposition~\ref{prop-Fell-multipliers} is
  satisfied, and condition~(ii) is not hard to check.  Condition~(iii)
  follows from the existence of an approximate identity for
  $C_{c}(\GG)$ in the inductive limit topology
  (cf. \cite[Corollary~2.11]{mrw}). Hence, $\pg$ extends to a
  nondegenerate homomorphism of $C^{*}(\GG)$ into
  $M(\cs(G\times\GG,\AA\times\GG))$ by
  Proposition~\ref{prop-Fell-multipliers}.

  Clearly, $\pa$ and $\pg$ commute.  Since $C^{*}(\GG)$
  is nuclear, we obtain a homomorphism $\pa\tensor\pg$ of
  $\cs(G,\AA)\tensor \cs(\GG)$ into $M(\cs(G\times\GG,\AA\times\GG))$.
   If $g\in\sa_c(G,\AA)$, $h\in C_c(\GG)$ and
    $k\in\sa_c(G\times\GG,\AA\times\GG)$, then an argument patterned
    after the proof of \cite[Lemma~1.108]{danacrossed} implies that
    $\pa(g)\pg(h)k$ is in $\sa_{c}(G\times\GG,\AA\times\GG)$ and that
    evaluation at $(s,x)\in G\times\GG$ ``passes through the
    integral'' in the second step in the next calculation:
    \begin{align*}
      \bigl( \pa(g)\pg(h)k\bigr)_1(s,x) &= \Bigl(\int_G
      \bigl(\pan(g(t))(\pg(h)k)\bigr)_1\,dt\Bigr)
      (s,x)\\
      &= \int_G \pan(g(t))\bigl(\pg(h)k\bigr)_1(s,x)\,dt\\
      &= \int_G g(t) \bigl(\pg(h)k\bigr)_1(t\inv s,x)\D(t)^{\half}\, dt\\
      &= \int_G\int_\GG g(t)h(y)k_1(t\inv s,y\inv
      x)d\lambda^{r(x)}(y)\D(t)^{\half}\, dt\\
      &= \int_{G\times\GG}\bigl((\D^{\half}g)\boxtimes h\bigr)_1(t,y)
      k_1((t,y)\inv(s,x))\,d\lambda^{r(s,x)}(t,y)\\
      &= \bigl( ((\D^{\half}g)\boxtimes h)*k\bigr)_1(s,x).
    \end{align*}%
  Therefore
  \begin{equation}\label{not hard}
    \pa\tensor\pg(g\tensor h)=(\D^{\half}g)\boxtimes h
    \quad\text{ for $g\in\sa_c(G,\AA)$ and $h\in C_c(\GG)$.}
  \end{equation}
  Since such elements $(\D^{\half}g)\boxtimes h$ span a dense subspace
  of $\sa_{c}(G\times\GG,\AA\times\GG)$, it follows that
  $\pa\tensor\pg$ maps $\cs(G,\AA)\tensor C^{*}(\GG)$ (into and) onto
  $\cs(G\times\GG,\AA\times\GG)$.

Now fix a faithful nondegenerate representation $\pi$ of
  $\cs(G,\AA)$ on a Hilbert space~$\HH$, and let $\pi_{0}:\AA\to
  B(\HH)$ be the nondegenerate representation whose integrated form
  is~$\pi$ (as in Lemma~\ref{bundle map}).  {Further let $\tau$
    be a faithful nondegenerate representation of $\cs(\GG)$ on a
    Hilbert space $\KK$. By the Disintegration Theorem
    (\cite[Proposition~4.2]{ren:representation} or
    \cite[Theorem~7.8]{muhwil:nyjm08}), we can assume
    $\KK=L^2(\GG^{(0)}*\VVV,\mu)$, where $\GG^{(0)}*\VVV$ is a Borel
    Hilbert bundle and $\mu$ is a finite quasi-invariant Radon measure
    on $\GG^{(0)}$, such that $\tau$ is the integrated form of a
    groupoid representation $\tau_0$ of $\GG$; thus
    \begin{equation}\label{tau}
      \bigl(\tau(h)\kappa\bigr)(u) = \int_{\GG} h(x)\tau_0(x)
      \kappa(s(x))\D_{\GG}(x)^{-\half}\, d\lambda^u(x) 
      \quad\text{ for $h\in C_c(\GG)$,}
    \end{equation}
    where $\D_{\GG}$ is the Radon-Nikodym derivative of $\nu\inv$ with
    respect to $\nu= \mu\circ\lambda$.  Note that we can identify
    $(G\times\GG)^{(0)}$ with $\go$.  Then we can form a Borel Hilbert
    bundle $\go*(\HH\tensor\VVV)$ such that
    $(\HH\tensor\VVV)(u)=\HH\tensor V(u)$ and such that
    $L^{2}(\go*(\HH\tensor\VVV),\mu)$ can be identified with $H\tensor
    L^{2}(\go*\VVV,\mu)$.  Then we can define a Borel $*$-functor (see
    \cite[Definition~4.5]{mw:fell}) $\Pi$ from $\AA\times\GG$ to
    $\operatorname{End}(\go*(\HH\tensor\VVV))$ by
    \begin{equation*}
      \Pi(a,x)=\pi_{0}(a)\tensor \tau_{0}(x).
    \end{equation*}
    If $\mu_{G}$ is a left Haar measure on $G$, then we get a Haar
    system $\set{\underline{\lambda}^{u}}_{u\in\go}$ on $G\times\GG$
    via $\underline{\lambda}^{u}=\mu_{G}\times\lambda^{u}$.  Notice
    that the Radon-Nikodym derivative of $\underline\nu^{-1}$ with
    respect to $\underline\nu:=\underline\lambda\circ\mu$ is given by
    $(s,x)\mapsto \Delta(s)\Delta_{\GG}(x)$.  Then
    \cite[Proposition~4.10]{mw:fell} implies that $\Pi$ integrates up
    to a $*$-homomorphism $L:\sa_{c}(G\times\GG,\AA\times\GG)\to
    B(\HH\tensor L^2(\GG^{(0)}*\VVV,\mu))$ given by
    \begin{multline}\label{L}
      L(f)(\eta\otimes\kappa)(u) \\ =\int_{G}\int_{\GG}
      \pi_{0}\bigl(f_{1}(t,x)\bigr)\eta\otimes\tau_0(x)\kappa(s(x))
      \D_{\GG}(x)^{-\half}\D(t)^{-\half}\,d\lambda^u(x)\, \,dt
    \end{multline}
    which extends to a representation of
    $\cs(G\times\GG,\AA\times\GG)$.

    Now, using~\eqref{not hard}, for $g\in \sa_c(G,\AA)$ and $h\in
    C_c(\GG)$ we have
    \begin{align*}
      L\bigl(&\rho_{\AA}\otimes\rho_{\GG}(g\otimes h)\bigr)
      (\eta\otimes\kappa)(u)\\
      &= L\bigl( (\D^{\half}g)\boxtimes h\bigr)(\eta\otimes\kappa)(u)\\
      &= \int_G \int_{\GG}
      \pi_0\bigl(\D^{\half}(t)g(t)h(x)\bigr)\eta\otimes\tau_0(x)\kappa(s(x))
      \D_{\GG}(x)^{-\half}\,d\lambda^u(x)\, \D(t)^{-\half}\,dt\\
      &= \Bigl( \int_G \pi_0(g(t))\eta\,dt\Bigr) \otimes\Bigl(
      \int_{\GG} h(x)\tau_0(x)
      \kappa(s(x))\D_{\GG}(x)^{-\half}\, d\lambda^u(x)\Bigr)\\
      &= \pi(g)\eta\otimes\bigl(\tau(h)\kappa\bigr)(u)\\
      &= (\pi\otimes\tau)(g\otimes h)(\eta\otimes\kappa)(u).
    \end{align*}
    It follows that $L\circ (\pa\tensor\pg)=\pi\tensor\tau$, and since
    the latter is a faithful representation of $\cs(G,\AA)\tensor
    C^{*}(\GG)$, it follows that $\pa\tensor\pg$ is faithful.}

  To complete the proof, we just let $\omega=(\pa\tensor\pg)^{-1}$.
  Then $\omega$ is an isomorphism of $\cs(G\times\GG,\AA\times\GG)$
  onto $\cs(G,\AA)\tensor C^{*}(\GG)$ and satisfies
  \begin{equation*}
    \omega(g\boxtimes h)=(\D^{-\half}g)\tensor h.\qed
  \end{equation*}
  \renewcommand\qed{\relax}
\end{proof}

\begin{proof}[Proof of Proposition~\ref{Upsilon}]
  Note that $C^{*}(\EE)=C^{*}(\EE,\lambda) \cong \KK(L^{2}(G))$. In
  fact, since $\EE$ is groupoid-equivalent to the trivial group,
  $C^{*}(\EE)$ is simple, so the representation $\tau:C^{*}(\EE)\to
  B(L^{2}(G))$ defined by
  \begin{equation*}
    \bigl(\tau(h)\kappa\bigr)(s) =\int_{G}h(s,t)\kappa(t)\,dt
    \quad\text{for $h\in C_{c}(\EE)$ and $\kappa\in
      C_{c}(G)\subseteq L^{2}(G)$}
  \end{equation*}
  is an isomorphism onto $\KK(L^{2}(G))$.  In particular, $C^*(\EE)$
  is nuclear, so by Lemma~\ref{omega}, we have an isomorphism
  \begin{equation*}
    \Upsilon :=(\id\tensor \tau)\circ
    \omega:\cs(G\times\EE,\AA\times\EE)\to \cs(G,\AA)\tensor 
    \KK\bigl(L^{2}(G)\bigr).
  \end{equation*}
  {If we let $\EE^{(0)}*\C$ be the trivial bundle
    $G\times\C$, then we can identify $L^{2}(G)$ with
    $L^2(\EE^{(0)}*\C,\lambda)$ in the obvious way.  Notice also that
    $\lambda$ is a quasi-invariant measure on $\EE^{(0)}$ with
    $\Delta_{\EE}\equiv 1$.  Thus the representation $\tau$ is
    essentially presented as in \eqref{tau}.  (The representation
    $\tau_{0}$ acts on $(s,z)\in G\times\C$ by
    $\tau_{0}(t,s)(s,z)=(t,z)$.)  Thus, in the current situation,
    \eqref{L} reduces to
    \[
    \bigl(L(f)\xi\bigr)(s) = \int_G \int_G
    \pi_0\bigl(f_1(r,s,t)\bigr)\xi(t)\D(r)^{-\half}\,dr\,dt
    \]
    for $\xi\in C_c(G,\HH)\subseteq \HH\otimes L^2(G)$.
    Now~\eqref{eq:9} is easily verified using the observation (from
    the proof of Lemma~\ref{omega}) that $(\pi\otimes\tau)\circ\omega
    = L$.}
\end{proof}

\begin{proof}[Proof of Theorem~\ref{Phi-iso}]
  We need to show that $\Phi$ is injective, and to do this we will
  show that the diagram
  \[
  \xymatrix{ \cs(G,\AA)\times_\d G\times_{\hat\d} G
    \ar[r]^-\Theta_-\cong \ar[d]_\P &\cs(\gltgg,\altgg)
    \ar[d]^\Psi_\cong
    \\
    \cs(G,\AA)\otimes \KK(L^2(G)) &\cs(\G\times\EE, \AA\times \EE)
    \ar[l]^-\Upsilon_-\cong }
  \]
  commutes, where $\Theta$, $\Psi$, and $\Upsilon$ are the
  isomorphisms of Propositions~\ref{Theta}, \ref{Psi},
  and~\ref{Upsilon}, respectively.

  Let $\pi:\cs(G,\AA)\to B(\HH)$ be a faithful nondegenerate
  representation on a Hilbert space $\HH$.  Let $f\in\sa_c(G,\AA)$,
  $g,h,\kappa\in C_c(G)$, and $\eta\in \HH$.  Then, to show that the
  diagram commutes, the following computation suffices.  Applying
  Proposition~\ref{Upsilon}, we have
  \begin{align*}
    \bigl((\pi\otimes{}&\id)\circ\Upsilon\circ\Psi\circ\Theta
    \bigl(k_\AA(f)k_{C(G)}(g)k_G(h)\bigr)(\eta\otimes\kappa)\bigr)(s)
    \\
    & =\int_G\int_G \pi_{0}\bigl(
    \Psi\circ\Theta\bigl(k_\AA(f)k_{C(G)}(g)k_G(h)\bigr)_1(r,s,t)
    \bigr)(\eta\otimes\kappa)(t)\D(r)^{-\half}\,dr\,dt
    \\
    \intertext{which, by Proposition~\ref{Psi}, is} & =\int_G\int_G
    \pi_{0}\bigl(
    \Theta\bigl(k_\AA(f)k_{C(G)}(g)k_G(h)\bigr)_1(r,r\inv s,s\inv rt)
    \bigr)\eta\kappa(t)\D(r)^{-\half}\,dr\,dt \\
    \intertext{which, by Proposition~\ref{Theta}, is} & =\int_G\int_G
    \pi_{0}\bigl( f(r)\D(r)^{\half}g(r\inv s)h(s\inv rt)\D(s\inv
    rt)^{\half}
    \bigr)\eta\kappa(t)\D(r)^{-\half}\,dr\,dt \\
    \intertext{which, after using Fubini and sending $t\mapsto
      r^{-1}st$, is} & =\int_G\int_G \pi(f(r))g(r\inv
    s)h(t)\D(t)^{\half} \eta\kappa(r\inv st)\,dt\,dr
    \\
    \intertext{which, since $\rho_t\kappa(r\inv s) = \kappa(r\inv
      st)\D(t)^{\half}$, is} & =\int_G\int_G \pi_{0}(f(r))\eta g(r\inv
    s) (\r_t\kappa)(r\inv
    s)h(t)\,dt\,dr \\
    & =\int_G \pi_{0}(f(r))\eta g(r\inv s)
    \bigl(\r(h)\kappa\bigr)(r\inv s)\,dr \\
    & =\int_G \pi_{0}(f(r))\eta
    \bigl(M_g\r(h)\kappa\bigr)(r\inv s)\,dr \\
    & =\int_G \pi_{0}(f(r))\eta \bigl(\l_rM_g\r(h)\kappa\bigr)(s)\,dr
    \\
    & =\int_G \bigl(\pi_{0}(f(r))\eta\otimes
    \l_rM_g\r(h)\kappa\bigr)(s)\,dr \\
    & =\int_G \bigl(\pi_{0}(f(r))\otimes
    \l_rM_g\r(h)\bigr)(\eta\otimes\kappa)(s)\,dr \\
    & =\int_G \bigl((\pi\otimes\l)(f(r)\otimes r) (1\otimes
    M_g\r(h)\bigr)(\eta\otimes\kappa)(s)\,dr
    \\
    & =(\pi\otimes\id)\Bigl(\int_G \bigl((\id\otimes\l)\circ\d(f(r))
    (1\otimes M_g\r(h)\bigr)(\eta\otimes\kappa)(s)\,dr\Bigr)
    \\
    & =(\pi\otimes\id)(\id\otimes\l)\circ\d(f)\bigl(1\otimes
    M_g\r(h)\bigr)(s)
    \\
    & =(\pi\otimes\id)\circ\P\bigl(k_\AA(f)k_{C(G)}(g)k_G(h)\bigr)(s).
    \qed
  \end{align*}
  \renewcommand\qed{\relax}
\end{proof}


\providecommand{\bysame}{\leavevmode\hbox to3em{\hrulefill}\thinspace}
\providecommand{\MR}{\relax\ifhmode\unskip\space\fi MR }
\providecommand{\MRhref}[2]{%
  \href{http://www.ams.org/mathscinet-getitem?mr=#1}{#2}
}
\providecommand{\href}[2]{#2}

\end{document}